\documentclass[microtype]{gtpart}    % Basic GT/GTM/AGT style
\usepackage{graphicx}  %%% the recommended graphics package
\usepackage[all]{xy}
\usepackage{enumerate}
\title[The algebraic duality resolution at $p=2$]{The algebraic duality resolution at $p=2$}

\author[A. Beaudry]{Agn\`es Beaudry}
\givenname{Agn\`es}
\surname{Beaudry}
\address{Department of Mathematics\\ University of Chicago \\ \newline 1118 East 58th Street \\ Chicago\\ Illinois \\ 60637}
\email{beaudry@math.uchicago.edu}
\urladdr{http://math.uchicago.edu/~beaudry/}


\arxivreference{math.AT/1412.2822}
\arxivpassword{6wbvm}

\volumenumber{}
\issuenumber{}
\publicationyear{}
\papernumber{}
\startpage{}
\endpage{}
\doi{}
\MR{}
\Zbl{}
\received{}
\revised{}
\accepted{}
\published{}
\publishedonline{}
\proposed{}
\seconded{}
\corresponding{}
\editor{}
\version{}

\newtheorem{thm}{Theorem}[subsection]
\newtheorem{cor}{Corollary}[subsection]
\newtheorem{prop}{Proposition}[subsection]
\newtheorem{lem}{Lemma}[subsection]
\theoremstyle{definition}
\newtheorem{defn}{Definition}[subsection]
\newtheorem{exmps}{Examples}[subsection]
\newtheorem{rem}{Remark}[subsection]

\makeatletter
\let\c@lem=\c@thm
\let\c@cor=\c@thm
\let\c@prop=\c@thm
\let\c@lem=\c@thm
\let\c@defn=\c@thm
\let\c@exmps=\c@thm
\let\c@rem=\c@thm
\makeatother

\numberwithin{equation}{subsection}
\newcommand{\bD}{\mathbb{D}}
\newcommand{\F}{\mathbb{F}}
\newcommand{\G}{\mathbb{G}}
\newcommand{\W}{\mathbb{W}}
\newcommand{\smsh}{\wedge}
\newcommand{\ra}{\rightarrow}
\newcommand{\xra}{\xrightarrow}

\DeclareFontFamily{OMS}{rsfs}{\skewchar\font'60}
\DeclareFontShape{OMS}{rsfs}{m}{n}{<-5>rsfs5 <5-7>rsfs7 <7->rsfs10 }{}
\DeclareSymbolFont{rsfs}{OMS}{rsfs}{m}{n}
\DeclareSymbolFontAlphabet{\scr}{rsfs}
\newcommand{\sB}{\scr{B}}
\newcommand{\sC}{\scr{C}}

\makeop{Gal}
\makeop{id}
\makeop{Mod}
\makeop{Tot}
\makeop{gr}
\makeop{Out}
\makeop{Hom}
\makeop{Ext}
\makeop{End}
\makeop{Aut}
\makeop{Tor}
\makeop{ev}
\begin{document}

\begin{abstract}    % type your abstract below
The goal of this paper is to develop some of the machinery necessary for doing  $K(2)$--local computations in the stable homotopy category using duality resolutions at the prime $p=2$. The Morava stabilizer group $\mathbb{S}_2$ admits a norm whose kernel we denote by $\mathbb{S}_2^1$. The algebraic duality resolution is a finite resolution of the trivial $\mathbb{Z}_2[[\mathbb{S}_2^1]]$--module $\mathbb{Z}_2$ by modules induced from representations of finite subgroups of $\mathbb{S}_2^1$. Its construction is due to Goerss, Henn, Mahowald and Rezk. It is an analogue of their finite resolution of the trivial $\mathbb{Z}_3[[\mathbb{G}_2^1]]$--module $\mathbb{Z}_3$ at the prime $p=3$. The construction was never published and it is the main result in this paper. In the process, we give a detailed description of the structure of Morava stabilizer group $\mathbb{S}_2$ at the prime $2$. We also describe the maps in the algebraic duality resolution with the precision necessary for explicit computations. 
\end{abstract}

\maketitle

%%%%%%%%%%%%%%%%%%%%   Start of main body of article

\section{Introduction}

Fix a prime $p$ and recall that $L_nS$ is the Bousfield localization of the sphere spectrum $S$ with respect to the wedge $K(0) \vee \ldots \vee K(n)$, where $K(m)$ is the $m$'th Morava $K$--theory at the prime $p$. The chromatic convergence theorem of Hopkins and Ravenel \cite[Section 8.6]{orangebook} states that the $p$--local sphere spectrum is the homotopy limit of its localizations $L_nS$. Further, there is a homotopy pull-back square:
\[\xymatrix{L_n S \ar[d] \ar[r]& L_{K(n)}S \ar[d] \\ L_{n-1}S \ar[r] &L_{n-1}L_{K(n)}S}\]
In theory, the homotopy groups of $S$ can be recovered from those of its Morava $K$--theory localizations $L_{K(n)}S$. For this reason, computing $\pi_*L_{K(n)}S$ is one of the central problems in stable homotopy theory. A detailed historical account of chromatic homotopy theory can be found in Goerss, Henn, Mahowald and Rezk \cite[Section 1]{ghmr}. 

The difficulty of computing $\pi_*L_{K(n)}S$ varies with $p$ and $n$. The computation of $\pi_*L_{K(1)}S$ is related to $K$--theory and is now well understood. For $n\geq 3$, almost nothing is known, which leaves the case $n=2$.  For $p\geq 5$, $\pi_*L_{K(2)}S$ was computed by Shimumora and Yabe in \cite{MR1318877}. Behrens gives an illuminating reconstruction of their results in \cite{MR2914955}. The case when $p=3$ proved much more difficult than the problem for $p\geq 5$. It is now largely understood due to the work of Shimomura, Wang, Goerss, Henn, Karamanov, Mahowald and Rezk (see Goerss and Henn \cite{GoerssBrown}, Goerss, Henn and Mahowald \cite{GoerssSplit}, Goerss, Henn, Mahowald and Rezk \cite{ghmr},  \cite{GoerssPic}, Henn, Karamanov and Mahowald \cite{HKM} and Shimomura and Wang \cite{MR1923218}).

The major breakthrough in understanding the case of $n=2$ and $p=3$ was the construction by Goerss, Henn, Mahowald and Rezk \cite{ghmr} of a finite resolution of the $K(2)$--local sphere called the \emph{duality resolution}. The duality resolution comes in two flavors. The \emph{algebraic} duality resolution is a finite resolution of the trivial $\Z_3[[\G_2]]$--module $\Z_3$ by permutation modules induced from representations of finite subgroups $G$ of the extended Morava stabilizer group $\G_2$. Its topological counterpart, the \emph{topological} duality resolution, is a finite resolution of $E_2^{h\G_2}$, where $E_2$ denotes Morava $E$--theory. It is composed of spectra of the form $\Sigma^k E_2^{hG}$, and realizes the algebraic duality resolution. Both the algebraic duality resolution and the topological duality resolution give rise to spectral sequences which can be used to study $\pi_*L_{K(2)}S$ at $p=3$. 

The existence of a resolution analogous to that of \cite[Theorem 4.1]{ghmr} at the prime $p=2$ was conjectured by Mahowald using the computations of Shimomura \cite{MR1702307} and of Shimomura and Wang \cite{MR1935487}. The central result of this paper is its construction, which is  due to Goerss, Henn, Mahowald and Rezk.  The author is grateful for their blessing to record it here. 

More precisely, for the norm one subgroup $\mathbb{S}_2^1$ (defined in (\ref{s21})), we construct a resolution of $\Z_2$ by modules which are induced from representations of finite subgroups of $\mathbb{S}_2^1$. We add a detailed description of the maps in the resolution, which will be used in later computations. However, we do not construct a full algebraic duality resolution of $\Z_2$ by $\Z_2[[\G_2]]$--modules as in \cite[Corollary 4.2]{ghmr} (see \fullref{rem:shortcomings}), nor do we realize the algebraic resolution topologically as in \cite[Section 5]{ghmr}. For work on the topological realization of the algebraic duality resolution, we refer the reader to Bobkova's thesis \cite{Bobkova}.

The algebraic duality resolution has already proved to be a useful tool for computations. We use the results of this paper in \cite{Paper2} to compute an associated graded for $H^*(\mathbb{S}_2^1, (E_2)_*V(0))$, where $V(0)$ is the mod $2$ Moore spectrum. The computations of \cite{Paper2} play a crucial role in \cite{Paper5}, where we prove that the strongest form of the chromatic splitting conjecture, as stated by Hovey in \cite[Conjecture 4.2 (v)]{cschov}, cannot hold when $n=p=2$.

\subsection{Background and notation}
As in Goerss, Henn, Mahowald and Rezk \cite[p.779]{ghmr}, ``we unapologetically focus on the case [$p = 2$] and $n = 2$ because this is at the edge of our current knowledge.''

We let $K(2)$ be Morava $K$--theory, so that
\[K(2)_* = \F_2[v_2^{\pm1}] \]
for $v_2$ of degree $6$, and whose formal group law is the Honda formal group law of height two, which we denote by $F_2$. The Morava stabilizer group $\mathbb{S}_2$ is the group of automorphisms of $F_2$ over $\F_{4}$. It admits an action of the Galois group $\Gal(\mathbb{F}_{4}/\mathbb{F}_2)$. The extended Morava stabilizer group $\G_2$ is
\[\G_2 = \mathbb{S}_2 \rtimes   \Gal(\mathbb{F}_{4}/\mathbb{F}_2).\]

By the Goerss--Hopkins--Miller theorem (see Goerss and Hopkins \cite{MR2125040}), the group $\G_2$ acts on Morava $E$--theory $E_2$ by maps of ${E}_{\infty}$--ring spectra and, for $X$ a finite spectrum, 
\begin{equation*} L_{K(2)}X \simeq E_2^{h\G_2} \smsh X. \end{equation*}
In fact, for any closed subgroup $G$ of $\G_2$, one can form the homotopy fixed point spectrum $E_2^{hG}$ (see the work of Hopkins and Devinatz  \cite{MR2030586} and of Davis \cite{MR2235364}). For a spectrum $X$, the action of $\G_{2}$ on $(E_{2})_*$ induces an action on
\begin{align*}
(E_{2})_*X := \pi_*L_{K({2})}(E_{2}\smsh X).
\end{align*}
For a closed subgroup $G$ of $\G_2$ and a finite spectrum $X$, there is a convergent descent spectral sequence
\begin{equation*}
E_2^{s,t} := H^s(G, (E_{2})_tX) \Longrightarrow \pi_{t-s}(E_2^{hG} \smsh X).
\end{equation*}
We describe the most relevant example for this paper here. There is a norm on the group $\mathbb{S}_2$ whose kernel is denoted $\mathbb{S}_2^1$ (see Goerss, Henn, Mahowald and Rezk \cite[Section 1.3]{ghmr}). Further,
\begin{align}\label{snsn1}
\mathbb{S}_2 \cong \mathbb{S}_2^1 \rtimes \Z_{2}.
\end{align}
Similarly, the norm on $\mathbb{S}_2$ induces a norm on $\G_2$. The kernel is denoted $\G_2^1$ and 
\begin{align}\label{gngn1}\G_2 \cong \G_2^1 \rtimes \Z_{2}.\end{align}
Let $\pi$ be a topological generator of $\Z_2$ in $\G_2$ and $\phi_{\pi}$ be its action on $E_2$. If $X$ is finite, there is a fiber sequence:
\begin{equation}\label{fseq}
L_{K(2)}X \ra E_2^{h\G_2^1}\smsh X \xra{\phi_{\pi}-\id} E_2^{h\G_2^1}\smsh X
\end{equation}
For this reason, the spectrum $E_2^{h\G_2^1}$ is often called the \emph{half sphere}. One approach for computing $\pi_*L_{K(2)}X$ is to compute the spectral sequence
\begin{equation}\label{dss1} H^s(\G_2^1, (E_2)_tX) \cong H^s(\mathbb{S}_2^1, (E_2)_tX)^{\Gal(\F_{4}/\F_2)}  \Longrightarrow \pi_{t-s}E_{2}^{h\G_2^1}\smsh X \end{equation}
and then use the fiber sequence (\ref{fseq}) to pass from $\pi_*(E_2^{h\G_2^1}\smsh X)$ to $\pi_*L_{K(2)}X$.

Computing the $E_2$--term of (\ref{dss1}) can be difficult. At the prime $3$, the algebraic duality resolution of Goerss, Henn, Mahowald and Rezk constructed in \cite[Theorem 4.1]{ghmr} gives rise to a first quadrant spectral sequence computing the $E_2$--term of the analogue of (\ref{dss1}). One of the most important consequences of this paper is the existence of such a spectral sequence at the prime $p=2$ (\fullref{thm:DRSS} below).

To state the main results and describe this spectral sequence, we must introduce some subgroups of $\mathbb{S}_2$. The group $\mathbb{S}_2$ has a unique conjugacy class of maximal finite subgroups of order $24$. Fix a representative and call it $G_{24}$. The group $G_{24}$ is isomorphic to the semi-direct product of a quaternion group denoted $Q_8$ and a cyclic group of order $3$ denoted $C_3$ (see \fullref{subsec:conj}). That is
\[G_{24} \cong Q_8 \rtimes C_3.\]
However, there are two conjugacy classes of maximal finite subgroups in $\mathbb{S}_2^1$. If $\pi$ is as above (\ref{fseq}), the groups $G_{24}$ and 
\[G_{24}' := \pi G_{24} \pi^{-1}\]
are representatives for the distinct conjugacy classes. The group $\mathbb{S}_2$ also contains a central subgroup $C_2$ of order $2$ generated by the automorphism $[-1](x)$ of the formal group law $F_2$ of $K(2)$. Therefore, $\mathbb{S}_2^1$ contains a cyclic subgroup 
\[C_6 := C_2 \times C_3 .\]
Choose a generator $\omega$ of $C_3$ and an element $i$ in $G_{24}$ such that $G_{24}$ is generated by $i$ and $\omega$. That is,
\[G_{24} = \left< i, \omega \right>.\]
Let $j = \omega i \omega^2$ and $k = \omega^2  j \omega$. The group $\mathbb{S}_2$ can be decomposed as a semi-direct product
\[ \mathbb{S}_2 \cong K \rtimes G_{24}\]
for a Poincar\'e duality group $K$ of dimension $4$.
Similarly, 
\[\mathbb{S}_2^1 \cong  K^1 \rtimes G_{24}\] 
for a Poincar\'e duality group $K^1$ of dimension $3$ (see \fullref{subsec:K}). The homology of the groups $K$ and $K^1$ play a central role in the construction of the duality resolution (see \fullref{subsec:theres}).

The group $\mathbb{S}_2^1$ is a profinite group and one can define the completed group ring 
\[ \Z_2[[\mathbb{S}_2^1]]  = \lim_{i,j} \Z/(2^i)[\mathbb{S}_2^1/U_j]\]
where $\{U_j\}$ forms a system of open subgroups such that $\bigcap_j U_j = \{e\}$. For any closed subgroup $G$ of $\mathbb{S}_2^1$, we let
\[\Z_2[[\mathbb{S}_2^1/G]] := \Z_2[[\mathbb{S}_2^1]] \otimes_{\Z_2[[G]]} \Z_2. \]

\subsection{Statement of the results}\label{subsec:state}
The main result of this paper is the following theorem.
\begin{thm}[Goerss--Henn--Mahowald--Rezk, unpublished]\label{thm:resintro}
Let $\Z_2$ be the trivial $\Z_2[[\mathbb{S}_2^1]]$--module. There is an exact sequence of complete left $\Z_2[[\mathbb{S}_2^1]]$--modules
\begin{align*}
0 \ra \sC_3  \xra{\partial_3} \sC_2  \xra{\partial_2}\sC_1  \xra{\partial_1}\sC_0 \ra \Z_2 \ra 0 
\end{align*}
where: 
\[\sC_p =  \begin{cases} \Z_2[[\mathbb{S}_2^1/G_{24}]] & p=0  \\ 
 \Z_2[[\mathbb{S}_2^1/C_{6}]]  & p=1,2\\
\Z_2[[\mathbb{S}_2^1/G_{24}']]  & p=3
\end{cases}\]
\end{thm}
The resolution of \fullref{thm:resintro} is called the \emph{algebraic duality resolution}. The name is motivated by the fact that the exact sequence of \fullref{thm:resintro} exhibits a certain twisted duality. To make this precise, let $\Mod({\mathbb{S}_2^1})$ denote the category of finitely generated left $\Z_2[[\mathbb{S}_2^1]]$--modules. As above, let $\pi$ be a topological generator of $\Z_2$ in $\mathbb{S}_2 \cong \mathbb{S}_2^1 \rtimes \Z_2$. For a module $M$ in $\Mod(\mathbb{S}_2^1)$, let $c_{\pi}(M)$ denote the left $\Z_2[[\mathbb{S}_2^1]]$--module whose underlying $\Z_2$--module is $M$, but whose $\Z_2[[\mathbb{S}_2^1]]$--module structure is twisted by the element $\pi$. 
\begin{thm}[Henn, Karamanov, Mahowald, unpublished]\label{thm:dualityintro}
Let 
\[\sC_p^*= \Hom_{\Z_2[[\mathbb{S}_2^1]]}(\sC_p, \Z_2[[\mathbb{S}_2^1]]).\] 
There is an isomorphism of complexes of left $\Z_2[[\mathbb{S}_2^1]]$--modules:
\begin{align*}
\xymatrix{\ar @{} [dr]  0 \ar[r] &  \sC_3 \ar[d]^{f_3} \ar[r]^-{\partial_3} & \sC_2 \ar[d]^{f_2} \ar[r]^-{\partial_2} & \sC_1  \ar[d]^{f_1} \ar[r]^-{\partial_1} & \sC_0 \ar[d]^{f_0} \ar[r] & \mathbb{Z}_2 \ar@{=}[d] \ar[r] & 0  \\ 0 \ar[r]	&  c_{\pi}(\sC_0^*) \ar[r]^-{c_{\pi}(\partial_1^*)}	& c_{\pi}(\sC_1^*) \ar[r]^-{c_{\pi}(\partial_2^*)}	& c_{\pi}(\sC_2^*) \ar[r]^-{c_{\pi}(\partial_3^*)}	& c_{\pi}(\sC_3^*) \ar[r]	& \mathbb{Z}_2 \ar[r]	 & 0	}
\end{align*}
\end{thm}

\begin{rem}\label{rem:shortcomings}
The resolution of \fullref{thm:resintro} has the following shortcoming: it does not extend to a resolution of the group $\mathbb{G}_2$ or of the group $\mathbb{S}_2$ as in \cite[Corollary 4.2]{ghmr}. This is due to the fact that (\ref{snsn1}) is a non-trivial extension when $n=p=2$. For $n=2$ and $p\geq 3$, $\mathbb{S}_2 \cong \mathbb{S}_2^1 \times \Z_p $.
\end{rem}

One application of the algebraic duality resolution is given by the following theorem.
\begin{thm}\label{thm:DRSS}
Let $M$ be a profinite left $\Z_2[[\mathbb{S}_2^1]]$--module. There is a first quadrant spectral sequence
 \[E_1^{p,q} = \Ext^q_{\Z_2[[\mathbb{S}_2^1]]}(\sC_p, M) \Longrightarrow H^{p+q}(\mathbb{S}_2^1,M)\]
with differentials $d_r \co E_r^{p,q}\ra E_r^{p+r,q-r+1}$.
Further, there are isomorphisms:
\[  E_1^{p,q}  \cong \begin{cases}
        H^q(G_{24}, M)  &p= 0\\
    H^q(C_6,M)   & p=1,2 \\
     H^q(G_{24}', M)  &p= 3\end{cases}  \]
    \end{thm}
The spectral sequence of \fullref{thm:DRSS} is called the algebraic duality resolution spectral sequence. Its computational appeal is twofold. The $E_1$--term is composed of the cohomology of finite groups. Further, it collapses at the $E_4$--term. The $d_1$ differentials are induced by the maps of the exact sequence in \fullref{thm:resintro}. In order to compute the spectral sequence, it is necessary to have a detailed description of these maps, which we do in \fullref{thm:maps} below. 

The following result introduces some important elements in $\mathbb{S}_2^1$.
\begin{thm}\label{thm:topgenint}
There is an element $\alpha$ in $K^1$ such that $\mathbb{S}_2$ is topologically generated by the elements $\pi$, $\alpha$, $i$ and $\omega$. The group $\mathbb{S}_2^1$ is topologically generated by the elements $\alpha$, $i$ and $\omega$.
\end{thm}
To state the next result, for any element $\tau$ in $G_{24}$, let
\[\alpha_{\tau} = \tau \alpha \tau^{-1}\alpha^{-1}.\]
Let $S_2^1$ be the $2$--Sylow subgroup of $\mathbb{S}_2^1$. The group $S_2^1$ admits a decreasing filtration, denoted $F_{n/2}S_2^1$ which will be defined in \fullref{subsec:filt}.

\begin{thm}\label{thm:maps}
Let $e$ be the canonical generator of $\Z_2[[\mathbb{S}_2^1]]$ and $e_p$ be the canonical generator of $\sC_p$. For a subgroup $G$ of $S_2^1$, let $IG$ be the kernel of the augmentation $\varepsilon \co \Z_2[[G]] \rightarrow \Z_2$. The maps $\partial_p \co \sC_p \ra \sC_{p-1}$ of \fullref{thm:resintro} can be chosen so that:
\begin{enumerate}[(a)]
\item $\partial_1(e_1) = (e-\alpha)e_0$
\item $\partial_2(e_2) = \Theta e_1$ for an element $\Theta$ in $\Z_2[[\mathbb{S}_2^1]]^{C_3}$ such that
\[\Theta\equiv  e+\alpha+i+j+k-\alpha_i-\alpha_j-\alpha_k  \mod \mathcal{J}, \]
where $\mathcal{J}$ is the left ideal
\[\mathcal{J} = (IF_{4/2}K^1, (IF_{3/2}K^1)(IS_2^1),(IK^1)^7, 2(IK^1)^3, 4 IK^1, 8). \]
\item there are isomorphisms of $\Z_2[[\mathbb{S}_2^1]]$--modules $g_p\co \sC_p  \ra \sC_p$ and differentials 
\[\partial_{p+1}' \co \sC_{p+1} \ra \sC_p\]
such that
\begin{align*}
\xymatrix{\ar @{} [dr]  0 \ar[r] &  \sC_3 \ar[d]^{g_3} \ar[r]^-{\partial_3} & \sC_2 \ar[d]^{g_2} \ar[r]^-{\partial_2} & \sC_1  \ar[d]^{g_1} \ar[r]^-{\partial_1} & \sC_0 \ar[d]^{g_0} \ar[r] & \mathbb{Z}_2 \ar@{=}[d] \ar[r] & 0  \\ 0 \ar[r]	&  \sC_3 \ar[r]^-{\partial_3'}	&  \sC_2 \ar[r]^-{\partial_2'}	&  \sC_1 \ar[r]^-{\partial_1'}		&\sC_0 \ar[r]& \mathbb{Z}_2 \ar[r]	 & 0	}
\end{align*}
is an isomorphism of complexes.
The map $\partial_3'\co \sC_3 \ra \sC_2$ is given by
\begin{equation*}
\partial_3'(e_3) = \pi (e+i+j+k)(e-\alpha^{-1}) \pi^{-1} e_2.
\end{equation*}
\end{enumerate}
\end{thm}
\fullref{thm:maps} is the key to doing computations using the duality resolution spectral sequence. The most difficult part of \fullref{thm:maps} is giving a good estimate for $\partial_2 \co \sC_2 \ra \sC_1$. A detailed description of the map $\partial_2$ is given in \fullref{subsec:maps}. Though such precision is not needed for our computations in \cite{Paper2}, the hope is that it will be sufficient for most future computations using the duality resolution spectral sequence.

\subsection{Organization of the paper}

\fullref{sec:mor} is dedicated to the description of the Morava stabilizer group in the special case of $p=2$ and $n=2$. 
(A more general account of the structure of $\mathbb{S}_n$ can be found in Goerss, Henn, Mahowald and Rezk \cite[Section 1]{ghmr}.) We begin by recalling the standard filtration on $\mathbb{S}_2$ and defining the norm. This allows us to define the unit norm subgroup $\mathbb{S}_2^1$ and describe its finite subgroups. In particular, we give an explicit choice of maximal finite subgroup $G_{24}$ in \fullref{lem:Q8split}. In \fullref{subsec:K}, we introduce a subgroup $K$ such that $\mathbb{S}_2 \cong K \rtimes G_{24}$ and compute the cohomology of $K$ and of its norm one subgroup $K^1$. These results are due to Goerss and Henn but are not published. We finish the section with a proof of \fullref{thm:topgenint}.

In \fullref{sec:algres}, we introduce the finite resolution of the trivial $\Z_2[[\mathbb{S}_2^1]]$--module $\Z_2$. We construct the algebraic duality resolution spectral sequence. We describe the duality properties of the resolution and give a proof of \fullref{thm:dualityintro}. We end this section by giving a detailed description of the maps in the resolution, proving \fullref{thm:maps}. 

The appendix, \fullref{sec:contcoh}, contains the results on the cohomology of profinite $p$--adic analytic groups used in this paper.

\subsection{Acknowledgements}
I thank Mark Behrens, Irina Bobkova, Daniel G. Davis, Matt Emerton, Niko Naumann, Doug Ravenel and Dylan Wilson for useful conversations. I thank the referee for a thorough reading and helpful comments. I thank Peter May for his support and for reading everything I have written multiple times. I thank Paul Goerss, Hans-Werner Henn and Charles Rezk for their generosity in sharing their ideas. In particular, I thank Paul and Hans-Werner for spending countless hours and hundreds of emails helping me understand their work and supporting mine. Finally, I will forever be indebted to Mark Mahowald for his incredible intuition and for the time he spent sharing it with me. I hope this paper and those to come will honor his tradition.

\section{The structure of the Morava stabilizer group}\label{sec:mor}

\subsection{A presentation of $\mathbb{S}_2$}\label{subsec:pres}
Let $F_2$ be the Honda formal group law of height $2$ at the prime $2$. It is the $2$--typical formal group law defined over $\F_2$ specified by the $2$--series
\[[2]_{F_2}(x) =x^4.\]
The ring of endomorphisms of $F_2$ over $\mathbb{F}_4$ is isomorphic to the maximal order $\mathcal{O}_2$ in the central division algebra over $\mathbb{Q}_2$ of valuation $1/2$, which we denote by
\[\mathbb{D}_2  = D(\mathbb{Q}_2, 1/2).\] 
We begin by describing this isomorphism. More details can be found in Ravenel \cite[A2.2]{greenbook} and Ravenel \cite[Chapter 4]{orangebook}.

Let $\W = W(\F_4)$ denote the ring of Witt vectors on $\F_4$.
The ring $\W$ is isomorphic to the ring of integers of the unique unramified degree $2$ extension of $\mathbb{Q}_2$. It is a complete local ring with residue field $\mathbb{F}_4$. The Teichm\"uller character defines a group homomorphism
\[ \tau \co \F_4^{\times} \ra \W^{\times}. \]
Let $\omega$ be a choice of primitive third root of unity in $\F_4^{\times}$, and identify $\omega$ with its Teichm\"uller lift $\tau(\omega)$. Given such a choice, there is an isomorphism
\[\W \cong \Z_2[\omega]/(1+\omega+\omega^2).\]
The Galois group $\Gal(\F_4/\F_2)$ is generated by the Frobenius $\sigma$. It is the $\Z_2$--linear automorphism of $\W$ determined by
\[\omega^{\sigma} = \omega^2.\]
The ring $\mathcal{O}_2$ is a non-commutative extension of $\W$
\[
\mathcal{O}_2 \cong  \W\left< S \right>/(S^2=2, aS=Sa^{\sigma})\]
for $a$ in $\W$. Note that any element of $\mathcal{O}_2$ can be expressed uniquely as a linear combination $a+bS$ for $a$ and $b$ in $\W$. The division algebra $\mathbb{D}_2$ is given by 
\[\mathbb{D}_2 \cong \mathcal{O}_2 \otimes_{\Z_2}\Q_2.\]
The $2$--adic valuation $v$ on $\Q_2$ extends uniquely to a valuation $v$ on $\bD_2$ such that $v(S) = \frac{1}{2}$. Further, $\mathcal{O}_2 = \{x \in \bD \mid v(x) \geq 0 \}$ and $ \mathcal{O}_2^{\times} =\{x \in \bD \mid v(x) =0 \} $. Therefore, any finite subgroup $G \subseteq \bD_2^{\times}$ is contained in $ \mathcal{O}_2^{\times}$.

Next, we describe the ring of endomorphisms of $F_2$ and give an explicit isomorphism $\End(F_2) \cong \mathcal{O}_2$. A complete proof can be found in Ravenel \cite[Section A2]{greenbook}. First, note that 
\[ \End(F_2) \subseteq \F_4[[x]].\]
To avoid confusion with the elements $\W \subseteq \mathcal{O}_2$, let $\zeta$ in $\F_4$ be a choice of primitive third root of unity in the field of coefficients $\F_4$. Let $S(x)$ correspond to the endomorphism
\begin{gather*}
S(x) =x^2\\
\tag*{\text{so that}}
 [2]_{F_2}(x) = x^4 =  S(S(x)) = S^2(x). \end{gather*}
Define $\omega^i(x) = \zeta^i x$ and $0(x)=0$. Given an element $a$ in $\W$, one can write it uniquely as $a = \sum_{i=0}^{\infty} a_{i} 2^i$ where $a_i$ in $\W$ satisfies the equation
\[x^4-x = 0.\]
That is, $a_i$ is in $\{0,1,\omega, \omega^2\}$.
Let $\gamma = a+bS$ be an element of $\mathcal{O}_2$. Let $a= \sum_{i\geq 0} a_{2i}2^i$ and $b = \sum_{i\geq 0} a_{2i+1}2^i$. Using the fact that $S^2=2$, the element $\gamma$ can be expressed uniquely as a power series
\begin{align}\label{gamexp}
\gamma= a_0 +2a_2+4a_4+\ldots + (a_1 +2a_3+4a_6+\ldots)S = \sum_{i \geq 0} a_iS^i.
\end{align}
One can show that
\[\gamma(x) = a_0(x) +_{F_2}  a_1 (x^2) +_{F_2} a_2 (x^4) +_{F_2} \ldots +_{F_2}  a_i (x^{2^i}) +_{F_2} \ldots  \]
is a well-defined power series in $\F_4[[x]]$. This determines a ring isomorphism from $\mathcal{O}_2$ to $\End(F_2)$.

The Morava stabilizer group $\mathbb{S}_2$ is the group of automorphisms of $F_2$. Thus, 
\[\mathbb{S}_2 \cong \mathcal{O}_2^{\times}.\] 
Any element of $\mathbb{S}_2$ can be expressed uniquely as a linear combination $a+bS$ for $a$ in $\W^{\times}$ and $b$ in $\W$. The center of $\mathbb{S}_2$ is given by the Galois invariant elements in $\W^{\times}$,
\[Z(\mathbb{S}_2) \cong \Z_2^{\times}.\]
Further, the element $\omega$ in $\W^{\times}$ generates a cyclic group of order $3$ in $\mathbb{S}_2$, denoted $C_3$. The reduction of $\W$ modulo $2$ induces an isomorphism $C_3 \cong \F_4^{\times}$.

The Galois group acts on $\mathbb{S}_2$ by 
\[(a+bS)^{\sigma} = a^{\sigma}+b^{\sigma}S.\]
The extended Morava stabilizer group $\G_2$ is defined by
 \[\G_2 := \mathbb{S}_2 \rtimes \Gal(\F_4/\F_2).\] 

\subsection{The filtration}\label{subsec:filt}
In what follows, we use the presentation of $\mathbb{S}_2$ induced by the isomorphism $\mathbb{S}_2 \cong \mathcal{O}_2^{\times}$ that was described in \fullref{subsec:pres}. That is,
\[ \mathbb{S}_2 \cong \left(\W\left< S \right>/(S^2=2, aS=Sa^{\sigma})\right)^{\times} \]
for $a$ in $\W$. As described in Henn \cite[Section 3]{MR1604171}, the group $\mathbb{S}_2$ admits the following filtration. 

Recall from \fullref{subsec:pres} that there is a valuation $v$ on $\mathcal{O}_2$ such that
\[v(S) = \frac{1}{2}.\]
Regard $\mathbb{S}_2$ as the units in $\mathcal{O}_2$. For all $n\geq 0$, define
\[F_{n/2}\mathbb{S}_2 = \{ x\in \mathbb{S}_2 \mid  v(x-1) \geq n/2\}.\]
This filtration corresponds to the filtration of $\mathbb{S}_2$ by powers of $S$, that is, for $n\geq 1$,
\begin{align}\label{filt}
F_{n/2}\mathbb{S}_2 = \{\gamma \in \mathbb{S}_2 \mid \gamma = 1 +a_nS^n+\ldots  \}.
\end{align}
The motivation for indexing the filtration by half integers is that the induced filtration on $\Z_2^{\times}\subseteq \mathbb{S}_2$ is the usual $2$--adic filtration by powers of $2$. 

Let
\begin{gather*}
\gr_{n/2}\mathbb{S}_2 := F_{n/2}\mathbb{S}_2 / F_{(n+1)/2}\mathbb{S}_2 \\
\tag*{\text{and}}
\gr\mathbb{S}_2 = \bigoplus_{n \geq 0}\gr_{n/2}\mathbb{S}_2. \end{gather*}
Define $S_2 := F_{1/2}\mathbb{S}_2$. The group $S_2$ is the $2$--Sylow subgroup of $\mathbb{S}_2$. The map $\mathbb{S}_2 \ra \F_4^{\times}$ which sends $\gamma$ to $a_0$ has kernel $S_2$. It induces an isomorphism
\begin{align*}
\gr_{0/2}\mathbb{S}_2 \cong \mathbb{F}_4^{\times}.
\end{align*}
Suppose that $n>0$ and that $\gamma$ is an element of $F_{n/2}\mathbb{S}_2$, so that
\begin{align*}
\gamma = 1 +a_{n} S^{n}+ \ldots
\end{align*}
for $a_i$ as in (\ref{gamexp}). Let $\overline{\gamma}$ denote the image of $\gamma$ in $\gr_{n/2}\mathbb{S}_2$. For $n\geq 1$, the map defined by $\overline{\gamma} \mapsto a_n$ gives a group isomorphism
\begin{align*}
\gr_{n/2}\mathbb{S}_2 \cong \F_4.
\end{align*}
It follows from these observations that the subgroups $F_{n/2}\mathbb{S}_2$ form a basis of open neighborhoods for the identity, so that $\mathbb{S}_2$ is a profinite topological group. 

Given any subgroup $G$ of $\mathbb{S}_2$, the filtration on $\mathbb{S}_2$ induces a filtration on $G$, defined by 
$F_{n/2}G = F_{n/2}\mathbb{S}_2 \cap G$. Let 
\begin{align}\label{grG}
\gr G = \bigoplus_{n \geq 0}\gr_{n/2}G
\end{align} 
be the associated graded for this filtration.

The associated graded $\gr S_2$ has the structure of a restricted Lie algebra. The bracket is induced by the commutator in $S_2$ and the restriction is induced by squaring. In \cite[Lemma 3.1.4]{MR1604171}, Henn gives an explicit description of the structure of this Lie algebra. We record this result in the case when $p=2$ and $n=2$.
\begin{lem}[Henn]\label{lem:cent}
For $n,m>0$, let $a \in  F_{n/2}S_2$ and $b \in F_{m/2}S_2$. Let $\overline{a}$ be the image of $a$ in $\gr_{n/2}S_2$ and $\overline{b}$ be the image of $b$ in $\gr_{m/2}S_2$. If $[a,b]$ denote the commutator $aba^{-1}b^{-1}$, then
\[
\overline{[a,b]} \equiv \overline{a}\overline{b}^{2^n} + \overline{a}^{2^m}\overline{b}  \ \in  \ \gr_{(n+m)/2}S_2.
\]
If $P(a) = a^2$, then
\[ \overline{P(a)} \equiv \begin{cases}
 \overline{a}^3 & \text{$\in  \ \gr_{2/2}S_2$ \ \ \ if $n = 1$} \\
         \overline{a}+\overline{a}^2 & \text{$\in \ \gr_{4/2}S_2$ \ \ \ if $n = 2$} \\
        \overline{a} & \text{$\in \ \gr_{n/2+1}S_2$ \ if $n>2$.}\end{cases}  \] 
 \end{lem}

 \subsection{The norm and the subgroups $\mathbb{S}_2^1$ and $\mathbb{G}_2^1$}\label{subsec:norm}

The group $\mathbb{S}_2 \cong \mathcal{O}_2^{\times}$ acts on $\mathcal{O}_2$ by right multiplication. This gives rise to a representation $\rho\co \mathbb{S}_2 \ra GL_2(\W)$:
\begin{align}\label{rep-rho} \rho(a+bS) =  \left( \begin{matrix}
a & b \\
2b^{\sigma} & a^{\sigma} \end{matrix} \right)
\end{align}
The restriction of the determinant to $\mathbb{S}_2$ is given by
\begin{align}\label{det}
\det(a+bS) = aa^{\sigma}-2bb^{\sigma}.
\end{align}
This defines a group homomorphism $\det \co \mathbb{S}_2 \ra \Z_2^{\times}$. 

\begin{lem}\label{lem:normsurj}
The determinant $\det \co \mathbb{S}_2 \ra \Z_2^{\times}$ is split surjective. 
\end{lem}
Before proving this lemma, we introduce elements of $\mathbb{S}_2$ that will play a key role in the remainder of this paper and in future computations. First, let
\begin{equation}\label{pi}
\pi:= 1+2\omega.
\end{equation}
By Hensel's lemma, $\mathbb{Z}_2$ contains two solutions of $f(x) = x^2+7$. One of them satisfies 
\[\sqrt{-7}\equiv 1+4 \mod 8.\]
This allows us to define
\begin{equation}\label{alpha}
\alpha := \frac{1-2\omega}{\sqrt{-7}}.
\end{equation}
Note that the elements $\pi$ and $\alpha$ are in $\W^{\times} \subseteq \mathbb{S}_2$.

\begin{proof}[Proof of \fullref{lem:normsurj}.]
The group $\Z_2^{\times}$ is topologically generated by $-1$ and $3$. It suffices to show that these values are in the image of the determinant. A direct computation shows that $\det(\pi) = 3$ and that $\det(\alpha)=-1$. Since $\alpha$ and $\pi$ commute, they define a splitting.  
\end{proof}

\begin{defn}
The norm $N \co \mathbb{S}_2 \ra  \Z_2^{\times}/\{\pm1\}$ is the composite:
\begin{align*}
\mathbb{S}_2 \xra{\det} \Z_2^{\times} \ra \Z_2^{\times}/\{\pm1\}
\end{align*}
\end{defn}
For a prime $p$ and an integer $i\geq 1$, let $U_i = \{x \in \Z_p^{\times} \mid x \equiv 1 \mod p^i\}$.
For $p=2$, there is a canonical identification
\begin{equation}\label{U2def} \Z_2^{\times} \cong \{\pm1\} \times U_2.\end{equation}
Therefore, the image of the norm is canonically isomorphic to the group $U_2$. Further, the group $U_2$ is non-canonically isomorphic to the additive group $\Z_2$.

The subgroup $\mathbb{S}_2^1$ is defined by the short exact sequence:
\begin{equation}\label{s21}
1 \ra \mathbb{S}_2^1 \ra \mathbb{S}_2 \xra{N} \Z_2^{\times}/\{\pm1\} \ra 1
\end{equation}
Any element $\gamma$ such that $N(\gamma)$ is a topological generator of $\Z_2^{\times}/\{\pm1\}$ determines a splitting. The element $\pi$ defined in (\ref{pi}) is an example. This gives a decomposition:
\begin{align}\label{isoK}
\mathbb{S}_2 \cong \mathbb{S}_2^1 \rtimes \Z_2^{\times}/\{\pm1\} \cong  \mathbb{S}_2^1 \rtimes \Z_2
\end{align}
Note that the group $\mathbb{S}_2^1$ is closed in $\mathbb{S}_2$ as it is the intersection of the finite index subgroups which are the kernels of the norm followed by the projections $U_2 \ra \Z/2^n$ for $n\geq 0$.

The norm $N$ extends to a homomorphism
\[N\co \G_2 \ra \Z_2^{\times}/\{\pm1\} \times \Gal(\F_4/\F_2) \ra \Z_2^{\times}/\{\pm1\},\]
where the second map is the projection. The subgroup $\G_2^1$ is the kernel of the extended norm and
\begin{align}\label{semG}
\G_2 \cong \G_2^1 \rtimes \Z_2^{\times}/\{\pm1\} \cong  \G_2^1 \rtimes \Z_2 .
\end{align}
We note that there is no splitting which is equivariant with respect to the action of the Galois group.

The filtration on $\mathbb{S}_2$ induces a filtration on $\mathbb{S}_2^1$ and 
\begin{align}\label{S21Syl}
S_2^1 := F_{1/2}\mathbb{S}_2^1
\end{align}
is the $2$--Sylow subgroup of $\mathbb{S}_2^1$. 

\begin{rem}
Note that for odd primes $p$, there is a canonical isomorphism 
\[\Z_p^{\times} \cong  C_{p-1} \times U_1,\] 
where $C_{p-1}$ is a cyclic group of order $p-1$. The exact sequence analogous to (\ref{s21}) is given by
\[1 \ra \mathbb{S}_2^1 \ra \mathbb{S}_2 \ra \Z_p^{\times}/C_{p-1} \ra 1.\]
Further, it has a central splitting. Therefore, when $p$ is odd, the Morava stabilizer group is a product
 \[\G_2 \cong \G_2^1\times \Z_p^{\times}/C_{p-1}  \cong \G_2^1 \times \Z_p.\]
There is no central splitting at the prime $p=2$ and the extensions (\ref{isoK}) and (\ref{semG}) are non-trivial. 
\end{rem}

We will need the following result in \fullref{subsec:K} to prove \fullref{thm:Fn2}.
\begin{lem}\label{lem:assS21}
For $n\geq 1$,
\begin{align*}
\gr_{n/2}{S}_2^1 =\begin{cases}
\F_2 & \mbox{if $n$ is even}\\
         \F_4 & \mbox{if $n$ is odd}.\end{cases}\end{align*} 
\end{lem}
\begin{proof}
Let $F_{0/2} \Z_2^{\times}  = \Z_2^{\times}$ and for $n\geq 2$ even, 
\begin{align*}
F_{n/2} \Z_2^{\times}= F_{(n-1)/2} \Z_2^{\times}  := U_{n/2} =  \{\gamma \ | \ \gamma \equiv 1 \mod 2^{n/2} \}.
\end{align*}
Let $\gamma$ be in $\mathbb{S}_2$. Let $n\geq 2$ be even and suppose that $\gamma$ has an expansion of the form
\[
\gamma = 1+a_{n-1}S^{n-1}+a_{n}S^{n} \mod S^{n+1}.
\]
By (\ref{det}),
\begin{align*}
\det(\gamma) &\equiv 1+2^{n/2}(a_n+a_n^{\sigma}) + a_{n-1}a_{n-1}^{\sigma} 2^{n-1} \mod 2^{n/2+1}
\end{align*}
which is in $F_{n/2}\Z_2^{\times}$. Therefore, the determinant preserves this filtration. In fact, it induces short exact sequences of $\F_2$--vector spaces:
\begin{align*}
0\ra \gr_{n/2}\mathbb{S}_2^1 \ra \gr_{n/2}\mathbb{S}_2 \ra \gr_{n/2}\Z_2^{\times}\ra 0
\end{align*}
The result then follows from the fact that
\[\gr_{n/2}\mathbb{Z}_2^{\times} =\begin{cases}
\F_2 & \mbox{if $n$ is even}\\
         0 & \mbox{if $n$ is odd}\end{cases}           \]
and $ \gr_{n/2}\mathbb{S}_2  \cong (\F_2)^2$ for $n\geq 1$.
\end{proof}

\subsection{Finite subgroups of $\mathbb{S}_{2}$}\label{subsec:conj}

In this section, we describe the finite subgroups of $\mathbb{S}_2$ that will be used in the construction of the resolution of \fullref{thm:resintro}. We also prove that there are two conjugacy classes of maximal finite subgroups in $\mathbb{S}_2^1$. This will be used in the proof of \fullref{thm:resintro}.

\fullref{prop:conj-iff} is a special case of Hewett \cite[Theorem 1.4]{hewett1}.
\begin{prop}[Hewett]\label{prop:conj-iff}
Any maximal finite non-abelian subgroup of $\mathbb{S}_2$ is isomorphic to a binary tetrahedral group
\[G_{24} \cong Q_8 \rtimes C_3.\]
Here, $Q_8$ is the quaternion group
\[Q_8 \cong \left<i,j\mid i^2=j^2, iji=j\right>\] 
and the action of $C_3$ permutes $i$, $j$ and $ij$.
\end{prop}

Our next goal is to prove that there are two conjugacy classes of maximal finite subgroups in $\mathbb{S}_2^1$. To do this, we will need some preliminary results. Note that the classification of conjugacy classes of maximal finite subgroups of $\mathbb{S}_2$ is addressed in Hewett \cite{hewett2} and in Bujard \cite{bujard}. According to Bujard \cite[Remark 1.36]{bujard}, Hewett's \cite[Theorem 5.3]{hewett2} is incorrect. However, \cite[Theorem 1.35]{bujard} in the case $n=p=2$ is also stated incorrectly. A correct statement can be found in \cite[Theorem 4.30]{bujard}. To avoid confusion, we restate the results we need.

 \begin{prop}[Bujard]\label{prop:conjclasses}
 There is a unique conjugacy class of groups isomorphic to $Q_{8}$, respectively $G_{24}$, in $\mathbb{S}_2$.
 \end{prop}
 \begin{proof}
For $Q_8$, this is \cite[Lemma 1.25]{bujard}. For $G_{24}$, this is \cite[Theorem 1.28]{bujard}.
 \end{proof}

It will be useful to have explicit choices of subgroups $Q_8$ and $G_{24}$. The proof of the following lemma is a direct computation.
\begin{lem}[Henn]\label{lem:Q8split}
Let 
\[i := \frac{1}{1+2\omega}(1-\alpha S).\]
Define $j = \omega i \omega^2$ and $k = \omega^2 i \omega = ij$. The elements $i$ and $j$ generate a quaternion subgroup of $\mathbb{S}_2$, denoted $Q_8$. The elements $i$ and $\omega$ generate a subgroup isomorphic to $G_{24}$. Further, in $\mathbb{D}_2$,
\[\omega = -\frac{1+i+j+k}{2}.\]
 \end{lem}
 
For $H$ a subgroup of $G$, let $N_{G}(H)$ be the normalizer of $H$ in $G$. Let $C_{G}(H)$ be the centralizer of $H$ in $G$. Note that the element $1+i$ in $\mathbb{D}_2^{\times}$ is in $N_{\mathbb{D}_2^{\times}}(Q_{8})$. Since the valuation $v(1+i) = \frac{1}{2}$, the restriction of the valuation to the normalizer is surjective. Therefore, there is an exact sequence:
 \begin{align*}
 &1 \ra N_{\mathbb{S}_2}(Q_8)  \ra N_{\mathbb{D}_2^{\times}}(Q_8) \ra \frac{1}{2}\Z \ra 0
 \end{align*}
Since $\mathbb{D}_2 \cong \Q_2(i,j)$, it follows by the Skolem--Noether theorem that $\Aut(Q_8)$ can be realized by inner conjugation in $\bD_2^{\times}$. There is an exact sequence:
 \begin{align*}
 &1 \ra C_{\mathbb{D}_2^{\times}}(Q_8)  \ra N_{\mathbb{D}_2^{\times}}(Q_8) \ra \Aut(Q_8) \ra 0
 \end{align*}
 The next proposition describes which of these automorphisms can be realized by conjugation in $\mathbb{S}_2$.
 \begin{prop}[Henn]\label{prop:autq8}
The subgroup of $\Aut(Q_8)$ that can be realized by conjugation by an element of $\mathbb{S}_2$ is isomorphic to the alternating group $A_4$. It is generated by conjugation by the elements $i$, $j$ and $\omega$.
\end{prop}
 \begin{proof}
The group $\Aut(Q_8)$ is isomorphic to the symmetric group $\mathfrak{S}_4$. One verifies by a direct computation that conjugation by $i$, $j$ and $\omega$ generates a subgroup of $\Aut(Q_8)$ isomorphic to $A_4$. Let $\Out_{\mathbb{S}_2}(Q_8)$ be the group of automorphisms of $Q_8$ that can be realized by conjugation in $\mathbb{S}_2$. Since $C_{\mathbb{D}_2^{\times}}(Q_{8}) \cong \Q_2^{\times}$ and $C_{\mathbb{S}_2}(Q_8) \cong \Z_2^{\times}$, there is a commutative diagram
\[ \xymatrix{
 \Z_2^{\times}  \ar[d] \ar[r] &  \Q_2^{\times} \ar[d]  \ar[r]^{v} & \Z \ar[d] \\ 
 N_{\mathbb{S}_2}(Q_8) \ar[d]  \ar[r] & N_{\mathbb{D}_2^{\times}}(Q_8)  \ar[d] \ar[r]^{v} & \frac{1}{2}\Z \ar[d]  \\
\Out_{\mathbb{S}_2}(Q_8)  \ar[r]   & \mathfrak{S}_4 \ar[r] & \Z/2 }\]
where the columns and rows are short exact. Therefore, $\Out_{\mathbb{S}_2}(Q_8)  \cong A_4$.
\end{proof}

\begin{lem}\label{lem:normq8}
Let $G_{24} = Q_8 \rtimes C_3$. The normalizer of $Q_8$ in $\mathbb{S}_2$ is given by
\[N_{\mathbb{S}_2}(Q_8) \cong U_2 \times G_{24}.\]
\end{lem}
\begin{proof}
By \fullref{prop:autq8}, there is a short exact sequence:
\[ 1 \ra C_{\mathbb{S}_2}(Q_8) \ra N_{\mathbb{S}_2}(Q_8) \ra  A_4 \ra 1\]
The centralizer is the subgroup $\Z_2^{\times} \cong C_2 \times U_2$ of $\mathbb{S}_2$. Since $G_{24}$ is defined by the extension
\[1 \ra C_2 \ra G_{24} \ra A_4 \ra 1\]
and the elements of $U_2$ are in the centralizer of $G_{24}$, it follows that $N_{\mathbb{S}_2}(Q_8)$ is isomorphic to $U_2 \times G_{24}$.
\end{proof}

Since the image of the norm is torsion free, any finite subgroup $G$ of $\mathbb{S}_2$ is contained in the kernel $\mathbb{S}_2^1$. Therefore, $\mathbb{S}_2^1$ has the same maximal finite subgroups as $\mathbb{S}_2$. However, there are more conjugacy classes in $\mathbb{S}_2^1$. 
\begin{prop}\label{prop:LpiG24}
There are two conjugacy classes of maximal finite subgroups in $\mathbb{S}_2^1$. One is the conjugacy class of $G_{24}$ defined in \fullref{lem:Q8split}. The other is $\xi G_{24} \xi^{-1}$, where $\xi$ is any element such that $N(\xi)$ is a topological generator of $U_2$. 
\end{prop}
\begin{proof}
Let $\Z_2^{\times} \subseteq \mathbb{S}_2$ be the center. Define
\[\mathbb{S}_2^0 := \mathbb{S}_2^1 \times U_2,\]
where $U_2$ is as in (\ref{U2def}).
The restriction of the determinant to $U_2$
surjects onto $(\Z_2^{\times})^2$. Therefore, there is an exact sequence
\[1 \ra \mathbb{S}_2^0 \ra \mathbb{S}_2 \ra \Z_2^{\times}/(\{\pm 1\},(\Z_2^{\times})^{2}) \ra 1 \]
and $\mathbb{S}_2/\mathbb{S}_2^0 \cong \Z/2$. If $N(\xi)$ is a topological generator for $ \Z_2^{\times}/\{\pm 1\}$, then $\xi$ is a representative for the non-trivial coset in $\mathbb{S}_2/\mathbb{S}_2^0$.

By \fullref{prop:conjclasses}, there is a unique conjugacy class of subgroups isomorphic to $G_{24}$ in $\mathbb{S}_2$. Since conjugation by elements of the center $\Z_2^{\times}$ is trivial, any two conjugacy classes in $\mathbb{S}_2^1$ differ by conjugation by an element of $\mathbb{S}_2/\mathbb{S}_2^0 \cong \Z/2$. Therefore, there are at most two conjugacy classes.

Next, we show that the conjugacy classes of $G_{24}$ and $\xi G_{24} \xi^{-1}$ are distinct in $\mathbb{S}_2^1$. Conjugation acts on the $2$--Sylow subgroups; hence, it suffices to prove the claim for the subgroup $Q_8$ of $G_{24}$. Suppose that there exists an element $\gamma$ in $\mathbb{S}_2^1$ such that
\[\xi Q_8 \xi^{-1} = \gamma Q_8 \gamma^{-1}. \]
This would imply that $\gamma^{-1}\xi$ is in $N_{\mathbb{S}_2}(Q_8)$. By \fullref{lem:normq8}, $\gamma^{-1}\xi $ is a product $z \tau$ for $z$ in $U_2$ and $\tau$ in $G_{24}$. This implies that $\xi = \gamma z \tau$. However, $\gamma z \tau$ is in $\mathbb{S}_2^0$. This is a contradiction, since the residue class of $\xi$ in $\mathbb{S}_2/\mathbb{S}_2^0 $ is non-trivial. Therefore, $G_{24}$ and $\xi G_{24} \xi^{-1}$ represent distinct conjugacy classes in $\mathbb{S}_2^1$.
\end{proof}

A choice for $\xi$ is the element $\pi$ defined in (\ref{pi}). For the remainder of this paper, $G_{24}'$ will denote
\begin{align} \label{piG24}
G_{24}' := \pi G_{24} \pi^{-1},
\end{align}
so that $G_{24}$ and $G_{24}'$ are representatives for the two conjugacy classes of maximal finite subgroups in $\mathbb{S}_2^1$.

\subsection{The Poincar\'e duality subgroups}\label{subsec:K}
In this section, we introduce the subgroups $K$ and $K^1$ and we describe their continuous cohomology rings $H^*(K, \F_2)$ and $H^*(K^1, \F_2)$ as $G_{24}$--modules. The author learned the results of this section from Paul Goerss and Hans-Werner Henn. We refer the reader to \fullref{sec:contcoh} for details on the cohomology of a profinite group. 

Let $K$ be the closure of the subgroup of $\mathbb{S}_2$ generated by $\alpha$ (as defined in (\ref{alpha})) and $F_{3/2}\mathbb{S}_2$. That is,
\[ K = \overline{\left<\alpha,F_{3/2}\mathbb{S}_2 \right>}.\]
\begin{prop}
The subgroup $K$ is normal in $\mathbb{S}_2$. Further,
${S}_2 \cong K \rtimes Q_{8}$
and
$\mathbb{S}_2 \cong K \rtimes G_{24}$.
\end{prop}
\begin{proof}
There is an isomorphism $\mathbb{S}_2 \cong S_2 \rtimes C_3$ and $\alpha$ commutes with the group $C_3$. Further, for any element $\gamma$ in $S_{2}$, it follows from \fullref{lem:cent} that  the commutator $[\gamma, \alpha]$ is in $F_{3/2}\mathbb{S}_2$. Since $\mathbb{S}_2 \cong S_2 \rtimes C_3$, and $F_{3/2}\mathbb{S}_2$ is normal, $K$ is also normal. The quotient $S_2/K$ is a group of order $8$ generated by the image of the elements $i$ and $j$ defined in \fullref{lem:Q8split}. The inclusion of $Q_8$ followed by the projection to $S_2/K$ is an isomorphism. This defines a splitting. Similarly, the group $\mathbb{S}_2/K$ is a group of order $24$ generated by the image of $\omega$ and $i$, and this defines a splitting.
\end{proof}

\begin{cor}
If $K^1$ is the kernel of the norm restricted to $K$, then $\mathbb{S}_2^1 \cong K^1 \rtimes G_{24}$. 
\end{cor}
\begin{proof}
The elements $\alpha$ and $\pi$ are in the group $K$ since $\alpha^{-1}\pi$ is in $F_{3/2}\mathbb{S}_2$. Therefore, the norm restricted to $K$ is surjective and $ \mathbb{S}_2^1/K^1 \cong  \mathbb{S}_2/K$.
\end{proof}

Our next goal is to compute the group cohomology of $K$ and $K^1$. We will need a few preliminary results.

 \begin{prop}\label{prop:prof2adic}
Any open subgroup of $S_2$ or of $S_2^1$ is a profinite $2$--adic analytic group.
 \end{prop}
\begin{proof}
According to Dixon, Du Sautoy, Mann and Segal \cite[Theorem 8.1]{dixonsegal}, a topological group is $2$--adic analytic if and only if it has an open subgroup which is a finitely generated powerful pro--$2$ group. (By \cite[Definition 3.1]{dixonsegal}, a pro--$2$ group $H$ is powerful if the quotient $H/\overline{H^4}$ is abelian, where $H^4 = \left<h^4 \mid h\in H\right>$.)
By \fullref{lem:cent}, if $n\geq 3$, then $F_{n/2}S_2$ is topologically generated by any finite set of elements that surjects onto $F_{n/2}S_2/F_{(n+2)/2}S_2$. Further, the image of $P^2 \co F_{n/2}S_2 \rightarrow F_{(n+4)/2}S_2$ is dense by \fullref{lem:cent}. If $n\geq 4$, then 
\[[F_{n/2}S_2 , F_{n/2}S_2 ] \subseteq F_{(2n)/2}S_2 \subseteq F_{(n+4)/2}S_2.\] 
This implies that $F_{n/2}S_2$ is powerful for $n\geq 4$. Since any open subgroup $G$ of ${S}_2$ contains $F_{n/2}S_2$ for some large $n$, it is a profinite $2$--adic analytic group. 

The proof for open subgroups of $S_2^1$ is similar, using $F_{n/2}S_2^1$ instead of $F_{n/2}S_2$.
\end{proof}

By \fullref{prop:prof2adic}, open subgroups of $S_2$ and $S_2^1$ are compact $2$--adic analytic groups. This motivates our use of the following definition, which can be found in Symonds and Weigel \cite[Section 4]{SymWeig}.
\begin{defn}\label{defn:PDgroup}
Let $G$ be a compact $p$--adic analytic group. Then $G$ is a Poincar\'e duality group of dimension $n$ if $G$ has cohomological dimension $n$ and
\[ H^s(G, \Z_p[[G]]) \cong \begin{cases} \Z_p &s=n \\
0 & s\neq n \end{cases}\]
as abelian groups. The right $\Z_p[[G]]$--module $H^n(G, \Z_p[[G]])$ is denoted $D_p(G)$ and called the \emph{compact dualizing module}. If the action of $\Z_p[[G]]$ on $D_p(G)$ is trivial, the group $G$ is called \emph{orientable}.
\end{defn}
\begin{rem}\label{rem:poindual}
For a Poincar\'e duality group $G$ of dimension $n$, one can show that $H_n(G, D_p(G))$ is isomorphic to $\Z_p$ (see Symonds and Weigel \cite[Theorem 4.4.3]{SymWeig}). Given a choice of generator $[G]$ for $H_n(G, D_p(G))$, the cap product induces a natural isomorphism 
\[H^{n-\ast}(G, -) \xrightarrow{\cap [G]} H_{\ast}(G, D_p(G) \otimes_{\Z_p} -).\]
The following observations are useful to compute $D_p(G)$. Let $\phi_{D_p(G)}\co G \to \Z_p^{\ast}$ be the representation associated to the action of $G$ on $D_p(G)$. Let $L(G)$ be the $\Q_p$-Lie algebra associated to $G$, as defined in Lazard \cite[Definition V.2.4.2.5]{lazard}. The right conjugation action of $G$ on itself induces a natural right action on $L(G)$, and thus a homomorphism $Ad \co G \to \Aut(L(G))$.
By \cite[Corollary 5.2.5]{SymWeig}, if $G$ is $p$-torsion free, 
\[\phi_{D_p(G)}(g) = \det(Ad(g)).\]
\end{rem}
\begin{prop}\label{prop:dualtrivial}
If an open subgroup $U$ of $\mathbb{S}_n$ is a Poincar\'e duality group, then it is orientable.
\end{prop}
\begin{proof}
This is the argument given by Strickland in the proof of \cite[Proposition 5]{strickland}. For any open subgroup $U$ of $\mathbb{S}_n$, $L(U)$ is isomorphic to the central division algebra $\mathbb{D}_n$ over $\mathbb{Q}_p$ of valuation $1/n$. For $g$ in $U$, the action $Ad(g)$ is given by conjugation in $\mathbb{D}_n$, which has determinant one.
\end{proof}

The next result relies on Lazard's theory of groups which are \emph{\'equi--$p$--valu\'e}. We refer the reader who is unfamiliar with the theory of Lazard to Huber, Kings and Naumann \cite[Section 2]{naumann} for an overview of the terminology.
\begin{thm}\label{thm:Fn2}
For $n\geq 3$, the group $F_{n/2}S_2$ is a Poincar\'e duality group of dimension $4$. The continuous group cohomology $H^*(F_{n/2}S_2 , \F_2)$ is the exterior algebra generated by
\[H^1(F_{n/2}S_2 , \F_2) \cong  \Hom_{\F_2}(\gr_{n/2}S_2 \oplus \gr_{(n+1)/2}S_2, \F_2) \cong \F_2^4.  \]
Similarly, $F_{n/2}S_2^1$ is a Poincar\'e duality group of dimension $3$ and $H^*(F_{n/2}S_2^1 , \F_2)$ is the exterior algebra generated by
\[H^1(F_{n/2}S_2^1 , \F_2) =   \Hom_{\F_2}(\gr_{n/2}S_2^1 \oplus \gr_{(n+1)/2}S_2^1, \F_2) \cong \F_2^3.  \]
\end{thm}
\begin{proof}
We define a filtration $w \co  F_{n/2}S_2 \ra \mathbb{R}_{+}^* \cup \{\infty\}$ in the sense of Lazard \cite[Definition II.1.1.1]{lazard}.
Let $w(1) = \infty$. For $k\geq 0$ and $x \in F_{(n+2k)/2}S_2 \backslash F_{(n+2k+2)/2}S_2$, let $w(x) =\frac{n+2k}{2}$.
With this filtration, $F_{n/2}S_2$ is \'equi--$p$--valu\'e of rank $4$ in the sense of Lazard \cite[V.2.2.7]{lazard}, with $\gr F_{n/2}S_2$ generated by 
\[F_{n/2}S_2/F_{(n+2)/2}S_2 \cong \gr_{n/2}S_2 \oplus \gr_{(n+1)/2}S_2.\] 
To verify that $w$ is a filtration and that $F_{n/2}S_2$ is \'equi--$p$--valu\'e with respect to $w$, one uses the formulas of \fullref{lem:cent}, noting that the squaring map 
\[P \co F_{(n+2k)/2}S_2/F_{(n+2k+2)/2}S_2 \ra F_{(n+2k+2)/2}S_2/F_{(n+2k+4)/2}S_2  \]
is an isomorphism if and only if $n\geq 3$. The result then follows from \cite[Proposition V.2.5.7.1]{lazard}, which states that $H^*(F_{n/2}S_2,\F_2)$ is an exterior algebra on the $\F_2$--linear dual of 
\[F_{n/2}S_2/P(F_{n/2}S_2) \cong F_{n/2}S_2/F_{(n+2)/2}S_2 \cong \F_4^2  \cong \F_2^4.  \]
According to Symonds and Weigel \cite[Theorem 5.1.5]{SymWeig}, this also implies that $F_{n/2}S_2$ is a Poincar\'e duality group of dimension $4$ (note that in \cite{SymWeig}, the authors imply in the third footnote that they use the terms \emph{uniformly powerful pro--$p$} and \emph{\'equi--$p$--valu\'e} interchangeably.)

To prove the second claim, we use the same filtration, $F_{n/2}S_2^1$. By \fullref{lem:assS21}, 
\[F_{n/2}S_2^1/F_{(n+2)/2}S_2^1  \cong \F_4 \oplus \F_2 \cong \F_2^3 .  \proved\]
\end{proof}

Recall that we use the convention that
\[\alpha_{\tau} = [\tau, \alpha] = \tau\alpha\tau^{-1} \alpha^{-1}.\]
The following congruences will be used in the computations of this section:
\begin{align*}
i &\equiv 1+S \mod S^2 & j &\equiv 1+\omega^2 S \mod S^2 \\
-1&\equiv 1+ S^2 \mod S^4 &  \alpha &\equiv 1+\omega S^2 \mod S^4  \\
\alpha_i &\equiv 1+S^3 \mod S^4  & \alpha_j &\equiv 1+\omega^2 S^3 \mod S^4\\
 \alpha^2 &\equiv 1+ S^4 \mod S^5 & \alpha\pi &\equiv  1+\omega S^4 \mod S^5
\end{align*}
They are obtained by a direct computation using the definitions of $\pi$, $\alpha$, $i$ and $j$, which were given in (\ref{pi}), (\ref{alpha}) and \fullref{lem:Q8split}.

\begin{defn}\label{defn:namealpha}
Let
\begin{align*}
\alpha_{0} &= \alpha  & \alpha_{1} &= \alpha_{i}  &  \alpha_{2} &= \alpha_j  & \alpha_{3} &= \alpha^2  & \alpha_{4} &= \alpha\pi
\end{align*}
and let $x_s$ in $\Hom_{\F_2}( \gr \mathbb{S}_2, \F_2)$ be the function dual to the image of $\alpha_s$ in $\gr \mathbb{S}_2$. The conjugation action of an element $\tau$ on an element $g$ is denoted by $\tau_*(g)$.
\end{defn}

\begin{rem}\label{rem:actformula}
The action by conjugation by $\tau$ can be computed using \fullref{lem:cent} and the formula
\[[\tau, \gamma]\gamma = \tau_*(\gamma).\] 
\end{rem}

\begin{cor}\label{cor:F32}
The continuous group cohomology $H^*(F_{3/2}S_2 , \F_2)$ is the exterior algebra generated by
\[H^1(F_{3/2}S_2 , \F_2) \cong  \F_2\{ x_1, x_2,x_3,x_4\}  \]
for $x_s$ as in \fullref{defn:namealpha}. The action of $\alpha$ on $H^1(F_{3/2}S_2 , \F_2)$ is trivial. 

Similarly, $H^*(F_{3/2}S_2^1 , \F_2)$ is the exterior algebra generated by
\[H^1(F_{3/2}S_2^1 , \F_2) \cong  \F_2\{ x_1, x_2,x_3\}  \]
with a trivial action by $\alpha$.
\end{cor}
\begin{proof}
By \fullref{thm:Fn2}, $H^*(F_{3/2}S_2,\F_2)$ is an exterior algebra generated by the $\F_2$--linear dual of $F_{3/2}S_2/F_{5/2}S_2$. This group is generated by the image of $\alpha_1$, $\alpha_2$, $\alpha_3$ and $\alpha_4$ of \fullref{defn:namealpha}. Therefore, $H^*(F_{3/2}S_2,\F_2)$ is the exterior algebra generated by $\F_2\{x_1, x_2,x_3,x_4\}$.
Since $\alpha$ is in $F_{2/2}S_2$, if $g$ is in $F_{3/2}S_2$, the commutator $[\alpha, g]$ is in $F_{5/2}S_2$. Using \fullref{rem:actformula}, we conclude that the action of $\alpha$ on $H^1(F_{3/2}S_2,\F_2)$ is trivial. 

The second claim follows in the same way from the fact that $F_{3/2}S_2^1/F_{5/2}S_2^1$ is generated by the image of $\alpha_1$, $\alpha_2$ and $\alpha_3$.
\end{proof}

\begin{lem}\label{lem:H1KZ2}
For $\alpha_i$ as defined in \fullref{defn:namealpha}, and $\overline{\alpha}_i$ its image in $H_1(K , \Z_2)$, there is an isomorphism
\[H_1(K , \Z_2) \cong   \Z/4\{\overline{\alpha}_0 \} \oplus \Z/2\{\overline{\alpha}_1, \overline{\alpha}_2\} \oplus \Z_2\{\overline{\alpha}_4 \},\]
where $2\overline{\alpha}_0$ is the image of $\alpha_3=\alpha^2$. Similarly,
\[H_1(K^1 , \Z_2) \cong   \Z/4\{\overline{\alpha}_0 \} \oplus \Z/2\{\overline{\alpha}_1, \overline{\alpha}_2\}. \]
The conjugation action of $Q_8$ on $K$ factors through the quotient of $Q_8$ by the central subgroup $C_2$. The induced action on $H_1(K,\Z_2)$ is trivial on $\overline{\alpha}_4$ and given by
\begin{align*}
i_*(\overline{\alpha}_0) &=  \overline{\alpha}_0 + \overline{\alpha}_1 & j_*(\overline{\alpha}_0) &=  \overline{\alpha}_0 + \overline{\alpha}_2\\
i_*(\overline{\alpha}_1) &=  \overline{\alpha}_1 & j_*(\overline{\alpha}_1) &=  \overline{\alpha}_1 + 2\overline{\alpha}_0\\
i_*(\overline{\alpha}_2) &=  \overline{\alpha}_2 + 2\overline{\alpha}_0 & j_*(\overline{\alpha}_2) &=  \overline{\alpha}_2.
\end{align*}
Hence, $H_1(K^1, \Z_2)$ is generated by the image of $\alpha$ as a $G_{24}$--module.
\end{lem}
\begin{proof}
First, we prove that the group $[K,K]$ is dense in $F_{5/2}S_2^1$. Note that $K$ is contained in $F_{2/2}S_2$. Let $a$ and $b$ be in $K$. For $\overline{a}$ and $\overline{b}$ as in \fullref{lem:cent}, 
\[\overline{[a,b]} = \overline{a} \overline{b}^4 + \overline{a}^4 \overline{b}  \in gr_{4/2}S_2. \]
Since $\overline{a}$ and $\overline{b}$ are in $\F_4$ and $x^4 =x$ for all $x$ in $\F_4$, this implies that $\overline{[a,b]}=0$ in $\gr_{4/2}S_2$. Therefore, $[a,b] \in F_{5/2}S_2$. 

Since the norm is multiplicative, the elements of $[K,K]$ have norm one. Hence, ${[K,K]}$ is contained in $F_{5/2}S_2^1$. Further, the map from $[K,K]$ to $F_{5/2}S_2^1/F_{7/2}S_2^1$ induced by the inclusion is surjective. Indeed, $F_{5/2}S_2^1/F_{7/2}S_2^1$ is generated by the image of the elements
$[\alpha, [i, \alpha]]$, $[\alpha, [j, \alpha]]$, and $[[i, \alpha], [j, \alpha]]$, both of which are in $[K^1,K^1]$.
By \fullref{cor:F32}, this implies that the composite
\[[K^1, K^1] \hookrightarrow [K,K] \rightarrow H_1(F_{5/2}S_2^1, \F_2)\]
is surjective. According to Behrens and Lawson \cite[Theorem 2.1]{MR2244259}, it then follows from results of Koch and Serre that $[K^1, K^1]$ and $[K,K]$ are dense in $F_{5/2}S_2^1$. Therefore, 
\[\overline{[K^1,K^1]} = \overline{[K,K]}=F_{5/2}S_2^1.\]
Hence, $H_1(K,\Z_2) \cong K/F_{5/2}S_2^1$ and $H_1(K^1,\Z_2) \cong K^1/F_{5/2}S_2^1$.
Since $\alpha$ and $\pi$ are in $K$, the norm $N \co  K \ra \Z_2^{\times}/\{\pm 1\}$ is split surjective. The image of $\alpha_4 = \alpha \pi $ generates $ \Z_2^{\times}/\{\pm 1\} \cong \Z_2$. Therefore, 
\[H_1(K,\Z_2) \cong H_1(K^1,\Z_2) \oplus \Z_2\{\overline{\alpha}_4\}.\]
Finally, $H_1(K^1,\Z_2)$ is generated by the image of $\alpha_0= \alpha$, $\alpha_1 = \alpha_i$ and $\alpha_2=\alpha_j$. Since $\alpha_i$ and $\alpha_j$ are in $F_{3/2}S_2$, $\alpha_i^2$ and $\alpha_j^2$ are in $F_{5/2}S_2=\overline{[K^1,K^1]}$, hence their images must have order $2$ in $K^1/\overline{[K^1,K^1]}$. Finally, $ \alpha^2 \equiv 1+S^4$ modulo $S^5$, so that the image of $\alpha$ has order $4$ in $K^1/\overline{[K^1,K^1]}$. We conclude that
\[H_1(K^1,\Z_2) \cong  \Z/4\{\overline{\alpha}_0\} \oplus \Z/2\{\overline{\alpha}_1, \overline{\alpha}_2\} .\]

The action of $Q_8$ by conjugation factors through $C_2=\{\pm 1\}$ since $C_2$ is in the center of $\mathbb{S}_2$. The action of the generators $i$ and $j$ is computed using \fullref{rem:actformula} and the following relations, which hold modulo $S^5$:
\begin{align*}
 [i,\alpha] &\equiv \alpha_i  &  [i, \alpha_i]  &\equiv 1  &    [i, \alpha_j] &\equiv \alpha^2   &  [i, \alpha\pi] &\equiv 1   \\
  [j,\alpha]  &\equiv \alpha_j   &   [j, \alpha_i] &\equiv \alpha^2 &  [j, \alpha_j]&\equiv 1 &   [j, \alpha\pi]&\equiv 1  .
\end{align*} 
These relations are obtained from \fullref{lem:cent}.
\end{proof}

\begin{cor}\label{cor:KP}
The group $K$ is an orientable Poincar\'e duality group of dimension $4$ and the group $K^1$ is an orientable Poincar\'e duality group of dimension $3$. 
\end{cor} 
\begin{proof}
It is a theorem of Serre \cite[Section 1]{MR0180619} that the cohomological dimension of a $p$--torsion free profinite group $G$ is equal to the cohomological dimension of any of its open subgroups.
The group $K$ is a $2$--group, and by \fullref{lem:cent}, the squaring operation $P$ on $K$ has a trivial kernel. Therefore, $K$ is torsion free. The group $F_{3/2}S_2$ is an open subgroup of $K$. Hence, the cohomological dimension of $K$ is equal to the cohomological dimension of $F_{3/2}S_2$, so that $K$ has cohomological dimension $4$. Similarly, $K^1$ has cohomological dimension $3$ since it contains $F_{3/2}S_2^1$ as an open subgroup. According to Symonds and Weigel \cite[Proposition 4.4.1]{SymWeig}, a profinite group $G$ of finite cohomological dimension is a Poincar\'e duality group if and only if it contains an open subgroup which is a Poincar\'e duality group. Therefore, both $K$ and $K^1$ are Poincar\'e duality groups.

Since $K$ is an open subgroup of $\mathbb{S}_2$, it follows from \fullref{prop:dualtrivial} that it is orientable. It remains to prove that $K^1$ is orientable. Let $\Z_2^{\times}$ be the center of $\mathbb{S}_2$. Let $U_2 \subseteq \Z_2^{\times}$ be as in (\ref{U2def}). The group $H = K^1 \times U_2$ is an open subgroup of $\mathbb{S}_2$ and hence $H$ is orientable. Further, $L(H)   \cong L(K^1) \oplus L(U_2)$ and the action of $H$ preserves the summands. Recall from \fullref{rem:poindual} that the representation $\phi_{D_2(H)}\co H \to \Z_2^{\ast}$ is given by the determinant of the adjoint action of $H$ on $L(H)$. For $g$ in $H$
\[\det(Ad(g)) = \det\left(Ad(g)|_{L(K^1)}\right) \det\left(Ad(g)|_{L(U_2)}\right).\]
Since $U_2$ is abelian, $\det\left(Ad(g)|_{L(U_2)}\right) =1$. It follows from the orientability of $H$ that $ \det\left(Ad(g)|_{L(K^1)}\right) =1$. In particular, this holds for any $g$ in $K^1$ and the representation $\phi_{D_2(K^1)}$ is trivial.
\end{proof}

\begin{thm}[Goerss--Henn, unpublished]\label{thm:HKF2} As an $\F_2$--algebra,
\[
H^*(K,\F_2) \cong \mathbb{F}_2[x_0, x_1, x_2, x_4]/ (x_0^2, x_1^2 + x_0x_1, x_2^2 + x_0x_2, x_4^2),
\]
where $x_s$ has degree one and is as in \fullref{defn:namealpha}. Further, 
\[H^*(K^1,\F_2) \cong  \mathbb{F}_2[x_0, x_1, x_2]/ (x_0^2, x_1^2 + x_0x_1, x_2^2 + x_0x_2).
\]
The conjugation action of $Q_8$ factors through $Q_8/C_2 \cong C_2 \times C_2$. It is trivial on $x_0$ and $x_4$. On $x_1$ and $x_2$, it is described by
\begin{align*}
i_* (x_1) &= x_0 +x_1 & j_*( x_1) &= x_1 \\ 
i_* (x_2) &= x_2  & j_* (x_2) &= x_0 + x_2
\end{align*}
so that the induced representation on $H^1(K^1, \F_2)$ is isomorphic to the augmentation ideal $I(Q_8/C_2)$, and $H^2(K^1,\F_2)$ is isomorphic to the co-augmentation ideal $I(Q_8/C_2)^*$.

\end{thm}
\begin{proof}
The spectral sequence for the group extension
\[1 \ra F_{3/2}\mathbb{S}_2 \ra K \ra \Z/2\{\overline{\alpha}_0\} \ra 1\]
has $E_2$--term given by
\[
\mathbb{F}_2[x_0] \otimes E(x_1, x_2, x_3,x_4).
\]
It follows from \fullref{lem:H1KZ2} and \fullref{lem:sq1} of the appendix that $x_0^2=0$. Since $x_3$ is the function dual to the image of $\alpha^2$ in $\gr \mathbb{S}_2$, we have that $d_2(x_3) = x_0^2$. Using the isomorphism $H^1(K,\F_2) \cong \Hom(K, \F_2)$ and \fullref{lem:H1KZ2},
one computes that
\[H^1(K,\F_2) \cong \F_2\{x_0,x_1,x_2,x_4\}.\]
Hence, $d_r(x_i)=0$ for $i \neq 3$. All other differentials are determined by these differentials and
\[E_{3}  \cong E_{\infty} \cong E(x_0, x_1, x_2, x_4).\]
Similarly, the $E_{2}$--term for the extension
\[1 \ra F_{3/2}\mathbb{S}_2^1 \ra K^1\ra  \Z/2\{\overline{\alpha}_0\} \ra 1\]
is given by
$\mathbb{F}_2[x_0] \otimes E(x_1, x_2, x_3)$, and $E_{3}  \cong E_{\infty} \cong E(x_0, x_1, x_2)$.

Now we determine the multiplicative extensions. First, note that it follows from \fullref{lem:sq1} that $x_4^2=0$ since $x_4$ is dual to a class that lifts to the free class $\overline{\alpha}_4$ in $H_1(K, \Z_2)$. Similarly, $x_1^2$ and $x_2^2$ are non-zero since they lift to $2$--torsion classes $\overline{\alpha}_1$ and $\overline{\alpha}_2$ in $H_1(K, \Z_2)$. Therefore, $x_1^2$ and $x_2^2$ are linear combinations of $x_0x_1$ and $x_0x_2$. We will show that $x_1^2 = x_0x_1$. The proof that $x_2^2=x_0x_2$ is similar.

Let $N$ be the closure of the normal subgroup of $K^1$ generated by $F_{6/2}K^1$, $\alpha^2$ and $\alpha_j$. That is,
\[N = \overline{\left<F_{6/2}K^1, \alpha^2, \alpha_j\right>}.\]
Since $[K^1, F_{3/2}K^1]\subseteq F_{6/2}K^1$, and $[\alpha, \alpha_j] = \alpha_j^2$, the group $K^1/N$ is a group of order $8$ generated by the image $a$ of $\alpha$ and the image $b$ of $\alpha_i$. The order of $a$ is $2$ and the order of $b$ is $4$. Further, since $[\alpha, \alpha_i] = \alpha_i^2$,
the group $K^1/N$ is isomorphic to the dihedral group $D_8$. Now, note that
\[H_1(D_8,\F_2) \cong \F_2\{a,b\}.\]
It is proved in Adem and Milgram \cite[Chapter IV, Theorem 2.7]{adem} that
\[H^*(D_8, \F_2) \cong \F_2[x,y,w]/(xy),\]
where $x$ is the function dual to $a$ and $y$ is the function dual to $a+b$. Changing the basis of $H_1(D_8, \F_2)$ from $\left<a,a+b\right>$ to $\left<a,b\right>$, sends the basis $\left<x,y\right>$ of $H^1(D_8,\F_2)$ to the basis $\left<y_0, y_1\right> = \left<x+y,y\right>$. We obtain the following presentation
\[H^*(K^1/N, \F_2) \cong  \F_2[y_0,y_1,w]/(y_1^2+y_0y_1).\] 
The projection induces a map
\[ f\co H^*(K^1/N, \F_2) \rightarrow H^*(K^1, \F_2)\]
with $f(y_0)=x_0$ and $f(y_1)=x_1$. Therefore, $x_1^2+ x_0x_1=0$ in $H^2(K^1/N, \F_2)$. 

The action of $Q_8$ follows from \fullref{lem:H1KZ2}. The isomorphism between the representation $H^1(K^1,\F_2)$ and the representation $I(Q_8/C_2)$ defined by
\[0 \ra I(Q_8/C_2) \ra \F_2[Q_8/C_2] \xra{\epsilon} \F_2 \ra 0 \]
is given by sending $x_0$ to the invariant $e+i+j+ij$, $x_1$ to $e+j$ and $x_2$ to $e+i$.  
\end{proof}

We end this section with a description of the integral homology of $K^1$ that will be used heavily in the proof of \fullref{thm:resintro}.
\begin{cor}[Goerss--Henn, unpublished]\label{cor:CohK}
The integral homology of $K^1$ is given by:
\[ H_n(K^1, \Z_2) = \begin{cases}
 \Z_2   & \mbox{$n=0,3$}\\
       \Z/4 \oplus (\Z/2)^2  & \mbox{$n=1$}\\
        0 & \mbox{$n=2$} \end{cases}  \] 
\end{cor}
\begin{proof}
The result for $n=1$ is \fullref{lem:H1KZ2}. The homology $H_*(K^1,\F_2)$ is dual to $H^*(K^1,\F_2)$, computed in \fullref{thm:HKF2}. The groups $H_n(K^1,\Z_2)$ for $n=2,3$ are computed from the long exact sequence associated to 
\[0 \ra \Z_2 \xra{2} \Z_2 \ra \F_2 \ra 0,\]
using the fact that $H_n(K^1,\F_2)$ and $H_1(K^1, \Z_2)$ are known.
\end{proof}

We finish this section by proving \fullref{thm:topgenint}.
\begin{proof}[Proof of \fullref{thm:topgenint}.]
Since $\mathbb{S}_2 \cong S_2 \rtimes C_3$ and $\omega$ generates $C_3$, it suffices to show that $S_2$ is generated by $\pi$, $\alpha$, $i$ and $j =\omega i \omega^{-1}$. Further, according to Behrens and Lawson \cite[Theorem 2.1]{MR2244259}, it suffices to prove that the inclusion $\left<\pi, \alpha, i, j\right> \ra S_2$ induces a surjective map 
\[H_1(S_2, \F_2) \cong S_2/\overline{S_2^*},\]
where $S_2^*$ is the group $S_2^2[S_2,S_2]$. The claim then follows from the isomorphism $S_2 \cong K \rtimes Q_{8}$, the surjectivity of the map 
\[\left<\pi, \alpha,\alpha_i, \alpha_j\right> \rightarrow K/\overline{K^*}\]
 and the fact that $i$ and $j$ generate $Q_8$. The argument for $\mathbb{S}_2^1$ is similar.
\end{proof}

\section{The algebraic duality resolution}\label{sec:algres}

This section is devoted to the construction of the algebraic duality resolution and the description of its properties. We refer the reader to the appendix (\fullref{sec:contcoh}) for background on the cohomology of profinite groups.

\subsection{The resolution}\label{subsec:theres} From now on, we fix $p=2$. The goal of this section is to prove \fullref{thm:resintro}, which was stated in \fullref{subsec:state}, and is restated as \fullref{thm:res} below. The proof is broken into a series of results given in \fullref{lem:N0}, \fullref{lem:N1}, \fullref{lem:N2} and \fullref{thm:C3}. All results in this section are due to Goerss, Henn, Mahowald and Rezk.

Let $G_{24}$ be the maximal finite subgroup of $\mathbb{S}_2$ defined in \fullref{lem:Q8split}. Recall that $G_{24}' = \pi G_{24} \pi^{-1}$ for $\pi  = 1+2\omega$ in $\mathbb{S}_2$. It was shown in \fullref{prop:LpiG24} that there are two conjugacy classes of maximal finite subgroups in $\mathbb{S}_2^1$, and that $G_{24}$ and $G_{24}'$ are representatives.
Recall that $C_2 = \{ \pm 1\}$ is the subgroup generated by $[-1](x)$ and $C_6 = C_2 \times C_3$. The group $K^1$ is the Poincar\'e duality subgroup of $\mathbb{S}_2^1$ which was defined in \fullref{subsec:K}. 

\begin{lem}\label{lem:N0}
Let $\sC_0 = \Z_2[[\mathbb{S}_2^1/G_{24}]]$ with canonical generator $e_0$. Let $\varepsilon \co \sC_0 \ra \Z_2$ be the augmentation 
\begin{align*}
\varepsilon \co \Z_2[[\mathbb{S}_2^1/G_{24}]]  \ra \Z_2\otimes_{\Z_2[[\mathbb{S}_2^1]]}\Z_2[[\mathbb{S}_2^1]] \otimes_{\Z_2[G_{24}]}\Z_2\cong \Z_2.
\end{align*}
Let $N_0$ be defined by the short exact sequence
\begin{equation}\label{resN0} 0 \ra N_0 \ra \sC_0 \xra{\varepsilon} \Z_2 \ra 0.\end{equation}
Then $N_0$ is the left $\Z_2[[\mathbb{S}_2^1]]$--submodule of $\sC_0$ generated by $(e-\alpha)e_0$, for $e$ in $\mathbb{S}_2^1$ the unit and $\alpha$ as defined in (\ref{alpha}).
\end{lem}
\begin{proof}
Since $\mathbb{S}_2^1 \cong K^1 \rtimes G_{24}$, $\sC_0 \cong   \Z_2[[K^1]]$ as a $\Z_2[[K^1]]$--module. Therefore, $N_0$ is isomorphic to the augmentation ideal $IK^1$. \fullref{lem:Ggen} of the appendix implies that $H_1(K^1, \Z_2)\cong H_0(K^1,N_0)$, where an isomorphism sends the image of $g$ in $K^1/\overline{[K^1,K^1]}$ to the image of $e-g$ in $IK^1/\overline{{IK^1}^2}$. It was shown in \fullref{lem:H1KZ2} that $K^1/\overline{[K^1,K^1]}$ is generated by $\alpha$ as a $G_{24}$--module. This implies that, as a $G_{24}$--module, $H_0(K^1,N_0)$ is generated by the image of $(e-\alpha)e_0$. Therefore, the map $F \co \Z_2[[\mathbb{S}_2^1]] \ra N_0$ defined by $F(\gamma) = \gamma (e-\alpha)e_0$ induces a surjective map
\[\F_2\otimes_{\Z_2[[K^1]]} F \co  \F_2 \otimes_{\Z_2[[K^1]]}\Z_2[[\mathbb{S}_2^1]]  \ra \F_2 \otimes_{\Z_2[[K^1]]} N_0.\]
By \fullref{lem:ghmr} of the appendix, $F$ itself is surjective, and $(e-\alpha)e_0$ generates $N_0$ as an $\Z_2[[\mathbb{S}_2^1]]$--module.
\end{proof}

\begin{lem}\label{lem:N1}
Let $N_0$ be as in \fullref{lem:N0}. Let $\sC_1 = \Z_2[[\mathbb{S}_2^1/C_6]]$ with canonical generator $e_1$. There is a map 
$\partial_1 \co \sC_1 \ra N_0$
defined by 
\begin{equation}\label{myp1}
\partial_1(\gamma e_1) = \gamma (e-\alpha) e_0
\end{equation}
for $\gamma$ in $\Z_2[[\mathbb{S}_2^1]]$. 
Further, let $N_1$ be defined by the short exact sequence
\begin{equation}\label{resN1} 0 \ra N_1 \ra \sC_1 \xra{\partial_1} N_0 \ra 0, \end{equation}
and let $\Theta_0$ in $\Z_2[[\mathbb{S}_2^1]]$ be any element such that $\Theta_0 e_1 $ is in the kernel of $\partial_1$ and
\[ \Theta_0 e_1 \equiv (3+i+j+k)e_1 \mod (4, IK^1).\]
Then $\Theta_0 e_1$ generates $N_1$ over $\mathbb{S}_2^1$.
\end{lem}
\begin{proof}
The element $\alpha$ satisfies $\tau\alpha = \alpha \tau$ for $\tau \in C_6$. Therefore, the map $\partial_1$ given by (\ref{myp1}) is well-defined.

Let $N_1$ be the kernel of $\partial_1$. Note that $\Z_2[[\mathbb{S}_2^1/C_6]] \cong \Z_2[[K^1]]^4$ as $\Z_2[[K^1]]$--modules,  generated by $e_1$, $ie_1$, $je_1$ and $k e_1$. Therefore, there is an isomorphism of $G_{24}$--modules
\[H_0(K^1, \sC_1) \cong \Z_2[G_{24}/C_6].\]
As $\Z_2[[K^1]]$--modules, $H_0(K^1, \sC_1)  \cong \Z_2^4$ generated by the image of the classes $e_1$, $ie_1$, $je_1$ and $k e_1$. Since $N_0\cong IK^1$, \fullref{lem:Ggen} of the appendix and \fullref{cor:CohK} imply that
\[H_1(K^1, N_0) \cong H_2(K^1, \mathbb{Z}_2) =0.\]
Therefore, the long exact sequence on cohomology gives rise to a short exact sequence
\begin{align*}
0\ra H_0(K^1, N_1) \ra H_0(K^1, \sC_1) \ra H_0(K^1, N_0) \ra 0.
\end{align*}
By \fullref{lem:Ggen} of the appendix and \fullref{lem:H1KZ2},
\[H_0(K^1, N_0) \cong H_1(K^1,\Z_2) \cong  \Z/4 \oplus (\Z/2)^2,\] 
which is all torsion. Thus, we can identify $H_0(K^1, N_1)$ with a free submodule of $H_0(K^1, \sC_1)$. Further, it must have rank $4$ over $\Z_2$. This can be made explicit as follows.

The map $H_0(K^1, \partial_1)$ sends the residue class of $\tau e_1$ to that of $\tau(e-\alpha)e_0$. For $\tau$ in $G_{24}$, $\tau^{-1}e_0 =e_0$, hence $\tau(e-\alpha)e_0= (e-\tau_*(\alpha))e_0$, where $\tau_*(\alpha) = \tau \alpha \tau^{-1}$. Again, using the boundary isomorphism $H_1(K^1, \Z_2) \cong H_0(K^1, N_0)$ of \fullref{lem:Ggen}, the formulas of \fullref{lem:H1KZ2} together with the fact that $k=ij$ can be used to compute:
\begin{align*}
\partial_1( e_1) &\equiv \overline{\alpha} & 
\partial_1(i e_1) &\equiv \overline{\alpha}+ \overline{\alpha_i} & 
\partial_1(j e_1) &\equiv \overline{\alpha}+ \overline{\alpha_j}  &
\partial_1(k e_1) &\equiv 3\overline{\alpha}+ \overline{\alpha_i}+ \overline{\alpha_j}
\end{align*}
Here, $\overline{a}$ is the image of $a$ in $H_1(K^1, \Z_2)$. As $\alpha$ generates a group isomorphic to $\Z/4$, and $\alpha_i$ and $\alpha_j$ both generate groups isomorphic to $\Z/2$, a set of $\Z_2$ generators for the kernel of $H_0(K^1, \partial_1)$ is given by the the elements:
\begin{align*}
f_1 &= -4e_1 &  
 f_2&= 2(i-e)e_1 &
f_3 &= 2(j-e)e_1  & 
f_4&=  (k-i-j-e)e_1
\end{align*}
Let
\begin{align*}
f=(3e+i+j+k)e_1 \in H_0(K^1, N_1).
\end{align*}
Then $f$ generates $H_0(K^1, N_1)\cong \Z_2[G_{24}/C_6] $ as a $G_{24}$--module. Indeed, using the fact that $G_{24}/C_6 \cong Q_8/C_2$, one computes: 
\begin{align*}
f_1 &=  1/3(i+j+k-5)f  & 
f_2&= if-f  & 
f_3 &= jf-f &
f_4&=  -k(f+f_1)
\end{align*}
(Note that $-\tau$ denotes $(-1)\cdot \tau$ for the coefficient $-1$ in $\Z_2$, as opposed to the generator of the central $C_2$ in $Q_8$.)

Next, we show that if
\[ f' \equiv f \mod (4, IK^1),\]
then $f'$ also generates $H_0(K^1, N_1)$ as a $G_{24}$--module. To do this, note that $\Z_2[Q_8/C_2]$ is a complete local ring with maximal ideal $\mathfrak{m} =(2, IQ_8/C_2)$. Hence, any element congruent to $1$ modulo $\mathfrak{m}$ is invertible. Therefore, if 
$f' = f +\epsilon f$ for $\epsilon$ in $\mathfrak{m}$,
then $f'$ is also a generator. However, for $a$ in $H_0(K_1, \sC_1) $,
\begin{align*}
4ae_1 &= a\frac{1}{3} ((e-i)+(e-j)+(e-k)+2e)f.
\end{align*}
Hence, $ a\frac{1}{3} ((e-i)+(e-j)+(e-k)+2e)$ is in $\mathfrak{m}$. Therefore, $4H_0(K_1, \sC_1)$ is contained in $\mathfrak{m}f$. 

Let $\Theta_0$ in $\Z_2[[\mathbb{S}_2^1]]$ be such that
\[ \Theta_0 e_1 \equiv (3+i+j+k)e_1 \mod (4, IK^1) .\]
Let $F \co \Z_2[[\mathbb{S}_2^1]] \ra N_1$ be the map defined by $F(\gamma) = \gamma \Theta_0 e_1$.
It induces a surjective map
\[\F_2\otimes_{\Z_2[[K^1]]} F\co  \F_2 \otimes_{\Z_2[[K^1]]}\Z_2[[\mathbb{S}_2^1]]  \ra \F_2 \otimes_{\Z_2[[K^1]]} N_1.\]
By \fullref{lem:ghmr} of the appendix, $F$ itself is surjective, and $\Theta_0 e_1$ generates $N_1$ as an $\Z_2[[\mathbb{S}_2^1]]$--module.
\end{proof}

Define $tr_{C_3} \co \Z_2[[\mathbb{S}_2^1]] \ra  \Z_2[[\mathbb{S}_2^1]]$ to be the $\Z_2$--linear map induced by
\begin{equation}\label{trace}
tr_{C_3}(g) = g+\omega g\omega^{-1} + \omega^{-1} g\omega \end{equation}
for $g$ in $\mathbb{S}_2^1$ and $\omega$ our chosen generator of $C_3$.

\begin{lem}\label{lem:N2}
Let $\sC_2 = \Z_2[[\mathbb{S}_2^1/C_6]]$ with canonical generator $e_2$. Let $\Theta$ in $\Z_2[[\mathbb{S}_2^1]]$ satisfy:
\begin{enumerate}
\item $\tau \Theta = \Theta \tau$ for $\tau$ in $C_6$
\item  $\Theta e_1$ is in the kernel of $\partial_1 \co \sC_1 \ra \sC_0$
\item $ \Theta e_1 \equiv (3+i+j+k)e_1$ modulo $(4, IK^1) $
\end{enumerate}
Then the map of $\Z_2[[\mathbb{S}_2^1]]$--modules $\partial_2 \co \sC_2 \ra \sC_1$ defined by
\begin{equation}\label{myp2}
\partial_2(\gamma e_2) = \gamma \Theta e_2\end{equation}
surjects onto $N_1 = \ker(\partial_1)$. Further, if $N_2$ is defined by the exact sequence
\begin{equation}\label{resN2}
0 \ra N_2 \ra \sC_2 \xra{\partial_2}  N_1 \ra 0,
\end{equation}
then $N_2 \cong \Z_2[[K^1]]$ as $\Z_2[[K^1]]$--modules.
\end{lem}

\begin{proof}
Choose an element $\Theta_0$ which generates $N_1$ as in \fullref{lem:N1}. Recall that
$C_6 \cong C_2 \times C_3$
and that $C_2$ is in the center of $\mathbb{S}_2$. Therefore, for $tr_{C_3}$ as defined by (\ref{trace}),
\[\Theta = \frac{1}{3}tr_{C_3}(\Theta_0)\]
satisfies properties (1), (2) and (3). The map $\partial_2$ given by (\ref{myp2}) is well-defined and surjects onto $N_1$ by \fullref{lem:ghmr}. 

Let $N_2 \subseteq \sC_2$ be the kernel of $\partial_2$ as in the statement of the result. The map $\partial_2$ induces an isomorphism $H_0(K^1, \sC_2) \cong H_0(K^1, N_1)$. Hence, for all $n$,
\[H_{n}(K^1,N_2) \cong H_{n+1}(K^1, N_1) \cong H_{n+2}(K^1, N_0) \cong H_{n+3}(K^1, \Z_2).\]
This implies:
\begin{align*}
H_n(K^1, N_2) \cong \begin{cases}
 \Z_2  & \mbox{$n=0$}\\
      0& \mbox{$n>0$} \end{cases} 
\end{align*}
Choose an element $e'$ in $N_2$ such that $e'$ reduces to a generator of $\Z_2$ in $H_0(K^1, N_2)$. Define $\phi \co \Z_2[[K^1]] \ra N_2$ 
by $\phi(k) =ke'$. Then $\Tor_{0}^{\Z_2[[K^1]]}(\F_2, \phi)$ is an isomorphism, and $\Tor_{1}^{\Z_2[[K^1]]}(\F_2, \phi)$ is surjective. By \fullref{lem:ghmr} of the appendix, $\phi$ is an isomorphism of $\Z_2[[K^1]]$--modules.
\end{proof}

Splicing the exact sequences (\ref{resN0}), (\ref{resN1}) and (\ref{resN2}) gives an exact sequence
\begin{equation}\label{Mbef}
0 \ra N_2 \ra \sC_2 \ra \sC_1 \ra \sC_0 \ra \Z_2 \ra 0
\end{equation}
which is a free resolution of $\Z_2$ as a trivial $\Z_2[[K^1]]$--module. The next goal is to find an isomorphism  $N_2 \cong \Z_2[[\mathbb{S}_2^1/G'_{24}]]$, where $G_{24}'=\pi G_{24}\pi^{-1}$ represents the other conjugacy class of maximal finite subgroups in $\mathbb{S}_2^1$. To prove this, we will need a few results. Before stating these, we introduce some notation.

Let $G $ be a subgroup of $\mathbb{S}_2$ which contains the central subgroup $C_2$. We define
\begin{equation*}
PG := G/C_{2}.\end{equation*} 
We let
\begin{align}\label{A4}
A_4 &:= PG_{24}\\
\label{A4'}   A_4' &:= PG_{24}'.\end{align}
The choice of notation is justified by the fact that both of these groups are isomorphic to the alternating group on four letters. Note also that, since $C_2$ is central, $PC_{6} \cong C_3$ and
$P\mathbb{S}_2^1 \cong K^1 \rtimes A_4$. Therefore, for any $G$ which contains $C_2$,
\[\Z_2[[\mathbb{S}_2^1/G]] \cong \Z_2[[P\mathbb{S}_2^1 / P{G}]]\]
as $\Z_2[[\mathbb{S}_2^1]]$--modules. To prove that $N_2 \cong \Z_2[[\mathbb{S}_2^1/G_{24}']]$, it will thus be sufficient to prove that 
\[N_2 \cong \Z_2[[P\mathbb{S}_2^1/A_4']]\] 
as $\Z_2[[P\mathbb{S}_2^1]]$--modules. 

We showed in \fullref{cor:KP} that $K^1$ is a Poincar\'e duality group (see \fullref{defn:PDgroup}). Further, there is an isomorphism of $\Z_2[[K^1]]$--modules
\begin{align*}
  \Z_2[[P\mathbb{S}_2^1 / A_{4}']] \cong \Z_2[[K^1]].
 \end{align*}
Hence:
\begin{align}\label{HnK}
H^n(K^1, \Z_2[[P\mathbb{S}_2^1 / A_{4}']]) \cong \begin{cases}
\Z_2  & \mbox{$n=3$}\\
  0 & \mbox{otherwise}\end{cases} 
  \end{align}

\begin{lem}
The inclusion $\iota \co K^1 \rightarrow P\mathbb{S}_2^1$ induces an isomorphism
\[\iota^* \co H^3(P\mathbb{S}_2^1,\Z_2[[P\mathbb{S}_2^1 / A_{4}']])\ra H^3(K^1, \Z_2[[P\mathbb{S}_2^1 / A_{4}']] ).\]
\end{lem}
\begin{proof}
The action of $A_4$ on $H^3(K^1, \Z_2[[P\mathbb{S}_2^1 / A_{4}']])$ is trivial. This follows from the fact that there are no non-trivial one-dimensional representations of $A_{4}$. Indeed, 
\begin{align*}
\Hom(A_{4}, Gl_1(\mathbb{Z}_2)) &= H^1(A_{4}, \Z_2^{\times})
\end{align*}
and $H^1(A_{4}, \Z_2^{\times})=0$. Since $P\mathbb{S}_2^1 \cong K^1 \rtimes A_4$, there is a spectral sequence
\[H^p(A_{4},H^q(K^1, \Z_2[[P\mathbb{S}_2^1 / A_{4}']] ) )\Longrightarrow H^{p+q}(P\mathbb{S}_2^1, \Z_2[[P\mathbb{S}_2^1 / A_{4}']]).\]
Because the action of $A_{4}$ on $H^3(K^1, \Z_2[[P\mathbb{S}_2^1 / A_{4}']])$ is trivial, (\ref{HnK}) implies that the edge homomorphism
\[ H^3(P\mathbb{S}_2^1, \Z_2[[P\mathbb{S}_2^1 / A_{4}']] )\rightarrow H^0(A_{4},H^3(K^1, \Z_2[[P\mathbb{S}_2^1 /A_{4}']] ) )\]
induced by the inclusion $\iota \co K^1 \rightarrow P\mathbb{S}_2^1$ is an isomorphism.
\end{proof}

\begin{lem}\label{lem:surjeta}
There are surjections
\begin{align*}
\eta &\co \Hom_{\Z_2[[P\mathbb{S}_2^1]]}(N_2, \Z_2[[P\mathbb{S}_2^1 / A_4']] ) \rightarrow H^{3}(P\mathbb{S}_2^1,  \Z_2[[P\mathbb{S}_2^1 /  A_4']]  )\\
\eta' &\co \Hom_{\Z_2[[K^1]]}(N_2, \Z_2[[P\mathbb{S}_2^1 /  A_4']]   ) \rightarrow H^3(K^1, \Z_2[[P\mathbb{S}_2^1/  A_4']] )
\end{align*}
making the following diagram commute
\begin{align}\label{edgecomm}
\xymatrix{\ar @{} [dr]    \Hom_{\Z_2[[P\mathbb{S}_2^1]]}(N_2, \ \Z_2[[P\mathbb{S}_2^1 /  A_4']]  )  \ar[d]_-{{\iota^*}}  \ar[r]^-{{\eta}} &H^3(P\mathbb{S}_2^1,  \Z_2[[P\mathbb{S}_2^1 /  A_4']]  ) \ar[d]^-{\iota^*}  \\  \Hom_{\Z_2[[K^1]]}(N_2,  \Z_2[[P\mathbb{S}_2^1 / A_4']]  ) \ar[r]^-{\eta'} &H^3(K^1,  \Z_2[[P\mathbb{S}_2^1 /  A_4']] ) 	}
\end{align}
where $\iota^*$ is the map induced by the inclusion $\iota \co K^1 \ra P\mathbb{S}_2^1$.
\end{lem}
\begin{proof}
Let $\sB_p = \sC_p$ for $0\leq p<3$ and $\sB_3 = N_2$. Resolving $\sB_p$ by projective $\Z_2[[P\mathbb{S}_2^1]]$--modules gives rise to spectral sequences 
\begin{gather*}
E_1^{p,q} \cong \Ext_{\Z_2[[P\mathbb{S}_2^1]]}^q(\sB_p,  \Z_2[[P\mathbb{S}_2^1 / A_4']]   ) \Longrightarrow H^{p+q}(P\mathbb{S}_2^1,  \Z_2[[P\mathbb{S}_2^1 /  A_4']] )\\
\tag*{\text{and}}
F_1^{p,q} \cong \Ext_{\Z_2[[K^1]]}^q(\sB_p,  \Z_2[[P\mathbb{S}_2^1 / A_4']]   ) \Longrightarrow H^{p+q}(K^1,  \Z_2[[P\mathbb{S}_2^1 /  A_4']] ).
\end{gather*}
These are first quadrant cohomology spectral sequence, with differentials
\begin{gather*}
d_r\co E_r^{p,q}\rightarrow E_r^{p+r, q-r+1}\\
\tag*{\text{and}}
d_r\co F_r^{p,q}\rightarrow F_r^{p+r, q-r+1}.\end{gather*}
Further, $\iota \co K^1 \ra P\mathbb{S}_2^1$ induces a map of spectral sequences
\[ \iota^*\co E_r^{p,q} \ra F_r^{p,q}. \]
Let $\eta$ be the edge homomorphism
\[\eta \co E_1^{3,0} \ra H^{3}(P\mathbb{S}_2^1,  \Z_2[[P\mathbb{S}_2^1 /  A_4']]  ) \]
and let $\eta'$ be the edge homomorphism 
\[\eta' \co F_1^{3,0} \ra H^{3}(K^1,  \Z_2[[P\mathbb{S}_2^1 /  A_4']]  ). \]

First, note that since the modules $\sB_p$ are projective $\Z_2[[K^1]]$--modules, $F_{r}^{p,q}$ collapses with $F_{\infty}^{p,q}=0$ for $q>0$ so that 
\[ F_{\infty}^{3,0} \ra H^3(K^1; \Z_2[[P\mathbb{S}_2^1/A_4']]) \]
is surjective. Hence, $\eta'$ is surjective.

In order to show that ${\eta}$ is surjective, it is sufficient to show that $E_1^{3-q,q} = 0$ for $q>0$. For $q=1$ and $q=2$, this follows from the fact that $\Z_2[[P\mathbb{S}_2^1/C_3]]$ is a projective $\Z_2[[P\mathbb{S}_2^1]]$--module. Hence, if $q>0$, then
\[ \Ext_{\Z_2[[P\mathbb{S}_2^1]]}^{q}(\Z_2[[P\mathbb{S}_2^1/C_3]],  \Z_2[[P\mathbb{S}_2^1 / A_4']] ) =0.\]
It remains to show that 
\[E_1^{0,3}  = \Ext_{\Z_2[[P\mathbb{S}_2^1]]}^{3}(\sB_0,  \Z_2[[P\mathbb{S}_2^1 / A_4']] ) \]
is zero, where $\sB_0 = \Z_2[[P\mathbb{S}_2^1/A_4]]$.

Let $V \cong C_2\times C_2$ be the $2$--Sylow subgroup of $A_4$. Then
\begin{align*}
E_1^{0,3}&=\Ext_{\Z_2[[P\mathbb{S}_2^1]]}^{3}(\sB_0,  \Z_2[[P\mathbb{S}_2^1 /A_4']]) \\
&\cong H^3(A_{4},   \Z_2[[P\mathbb{S}_2^1 / A_4']] )\\
&\cong H^3(V,   \Z_2[[P\mathbb{S}_2^1 /A_4']] )^{C_3}.
\end{align*}

Let $G_n = PF_{n/2}S_2^1 \rtimes A_4'$ and $X_n = P\mathbb{S}_2^1 / G_n$. The profinite $A_4$--set $P\mathbb{S}_2^1 / A_4'$ is isomorphic to the inverse limit of the finite $A_4$--sets $X_n$. There is an exact sequence
\[0\ra \lim{}^1H^2(V,   \Z_2[X_n] ) \ra  H^3(V,   \Z_2[[P\mathbb{S}_2^1 / A_4']] ) \ra \lim_n H^3(V,   \Z_2[X_n]) \ra 0 . \]
Since the groups $H^2(V,   \Z_2[X_n] )$ are finite, the Mittag-Leffler condition is satisfied and $\lim{}^1H^2(V,   \Z_2[X_n] ) =0$.  Hence, 
\[H^3(V,   \Z_2[[P\mathbb{S}_2^1 / A_4']] )   \cong  \lim_n H^3(V,   \Z_2[X_n]).\]
We will show that there is an integer $N$ such that $H^3(V,   \Z_2[X_n]) =0$ for all $n \geq N$. This implies that $H^3(V,   \Z_2[[P\mathbb{S}_2^1 / A_4']] )$ is zero, so that $E_1^{0,3}=0$.

Note that
\begin{align*}
 H^3(V,   \Z_2[X_n]) &\cong \bigoplus_{x \in V \backslash X_n/G_n}H^3(V, \Z_2[V/V_x]) \\
 &\cong \bigoplus_{x \in V \backslash X_n /G_n}H^3(V_x, \Z_2)
\end{align*} 
for $V_{x} = \{g \in V \mid gxG_n = xG_n\}$. If the inclusion $V_{x} \subseteq V$ is an equality, then $x^{-1}Vx$ is a subgroup of $G_n$. We show that there exists an integer $N$ such that, for all $n \geq N$, there is no element $x$ in $P\mathbb{S}_2^1$ such that $x^{-1}Vx \subseteq G_n$. This implies that, for $n\geq N$, for all choices of coset representatives $x \in V \backslash X_n /G_n$, the group $V_x$ is either trivial or it has order $2$. In both cases, $H^3(V_x, \Z_2) =0$. Hence, for $n \geq N$, $H^3(V, \Z_2[X_n]) =0$.

Suppose that that there is a sequence of integers $n_m$ and elements $x_{n_m}$ such that $x_{n_m}^{-1}V{x_{n_m}} \subseteq G_{n_m}$. Since $P\mathbb{S}_2^1$ is compact, we can choose the sequence $(x_{n_m})$ to converge to some element $y$. The groups $G_n$ are closed and nested, so the continuity of the group multiplication implies that $y^{-1}Vy \subseteq G_n$ for all $n \in \mathbb{N}$. Therefore,
\[ y^{-1}Vy \subseteq \bigcap_n G_n = A_4'\]
and hence $y^{-1}Vy = V'$, where $V'$ is the $2$-Sylow subgroup of $A_4'$. However, it follows from \fullref{prop:LpiG24} that $V$ and $V'$ are not conjugate in $P\mathbb{S}_2^1$. Therefore, such a sequence cannot exist, and there must be some integer $N$ such that, for all $n \geq N$, there is no $x$ in $P\mathbb{S}_2^1$ such that $x^{-1}Vx \subseteq G_n$.  \end{proof}

\begin{thm}\label{thm:C3}
There is an isomorphism of left $\Z_2[[\mathbb{S}_2^1]]$--modules
\[\phi \co \Z_2[[\mathbb{S}_2^1/G_{24}']] \ra N_2,  \]
where $G_{24}' = \pi G_{24} \pi^{-1}$.
\end{thm}
\begin{proof}It suffices to construct an isomorphism $\varphi\co N_2 \to \Z_2[[P\mathbb{S}_2^1/PG_{24}']]$ of left $\Z_2[[P\mathbb{S}_2^1]]$--modules. The result then follows by letting $\phi = \varphi^{-1}$, considered as a map of $\Z_2[[\mathbb{S}_2^1]]$--modules.

Recall from \fullref{cor:KP} that $K^1$ is an orientable Poincar\'e duality group of dimension $3$, as in \fullref{defn:PDgroup}. That is, the compact dualizing module $D_2(K^1)$ is isomorphic to the trivial $\Z_2[[K^1]]$--module $\Z_2$ and $H_3(K^1,\Z_2) \cong \Z_2$. Choose a generator $[K^1]$ of $H_3(K^1,\Z_2)$. As in \fullref{rem:poindual}, there is a natural isomorphism 
\[H^{3-\ast}(K^1, -) \xrightarrow{\cap [K^1]} H_{\ast}(K^1, -).\]
Let
\[\nu \co H_3(K^1,\Z_2) \to \Z_2 \otimes_{\Z_2[[K^1]]}N_2\]
be the edge homomorphism for the homology spectral sequence obtained from (\ref{Mbef}). Define
\[\ev \co   \Hom_{\Z_2}(\Z_2 \otimes_{\Z_2[[K^1]]}N_2,   \Z_2 \otimes_{\Z_2[[K^1]]} \Z_2[[P\mathbb{S}_2^1 / A_4']] )  \to H_0(K^1,  \Z_2[[P\mathbb{S}_2^1 / A_4']]) \]
by
\[\ev(f) = f(\nu([K^1])).\]
Let $\iota \co K^1 \to P\mathbb{S}_2^1$ be the inclusion. Let $\eta$ and $\eta'$ be the edge homomorphisms of \fullref{lem:surjeta}. We obtain the following commutative diagram:
\[\xymatrix{ \Hom_{\Z_2[[P\mathbb{S}_2^1]]}(N_2,   \Z_2[[P\mathbb{S}_2^1 / A_4']]   )  \ar[r]^-{\eta} \ar[d]_-{\iota^*} & H^3(P\mathbb{S}_2^1,  \Z_2[[P\mathbb{S}_2^1 / A_4']] ) \ar[d]^-{\iota^*}   \\
 \Hom_{\Z_2[[K^1]]}(N_2,   \Z_2[[P\mathbb{S}_2^1 / A_4']] )  \ar[r]^-{\eta'} \ar[d]_-{\Z_2\otimes_{\Z_2[[K^1]]} -} & H^3(K^1,  \Z_2[[P\mathbb{S}_2^1 / A_4']] ) \ar[d]^-{\cap [K^1]}   \\
 \Hom_{\Z_2}(\Z_2 \otimes_{\Z_2[[K^1]]}N_2,   \Z_2 \otimes_{\Z_2[[K^1]]} \Z_2[[P\mathbb{S}_2^1 / A_4']] ) \ar[r]^-{\ev}   & H_0(K^1, \Z_2[[P\mathbb{S}_2^1 / A_4']]  )   } \]
Since $\cap [K^1] \circ \eta'$ is surjective, so is the map $\ev$. Both $N_2$ and $ \Z_2[[P\mathbb{S}_2^1 / A_4']]$ are free of rank one over $\Z_2[[K^1]]$. Hence, $\Z_2 \otimes_{\Z_2[[K^1]]}N_2$ and $ \Z_2 \otimes_{\Z_2[[K^1]]} \Z_2[[P\mathbb{S}_2^1 / A_4']] $ are abstractly isomorphic to $\Z_2$. Since $\ev$ is a surjective group homomorphism from $\Z_2$ to itself, it is an isomorphism. It follows from \fullref{lem:ghmr} that any element of $ \Hom_{\Z_2[[K^1]]}(N_2,   \Z_2[[P\mathbb{S}_2^1 / A_4']] )$ that becomes a unit after applying $\Z_2\otimes_{\Z_2[[K^1]]}-$ is an isomorphism. By \fullref{lem:surjeta}, the composite $\cap [K^1] \circ \iota^* \circ \eta$ is surjective. Therefore, we can choose an element $\varphi$ in $\Hom_{\Z_2[[P\mathbb{S}_2^1]]}(N_2,   \Z_2[[P\mathbb{S}_2^1 /  A_4']] )$ 
such that $\cap [K^1] \circ \iota^*\circ\eta(\varphi)$ is a generator of $H_0(K^1, \Z_2[[P\mathbb{S}_2^1 /  A_4']] )$. Then $\iota^*(\varphi)$ in $\Hom_{\Z_2[[K^1]]}(N_2,   \Z_2[[P\mathbb{S}_2^1 /  A_4']] )$ is an isomorphism, and hence $\varphi$ must be an isomorphism.
\end{proof}

Combining the previous results, we can finally prove \fullref{thm:resintro}. We restate it here for convenience. 

\begin{thm}\label{thm:res}
Let $\Z_2$ be the trivial $\Z_2[[\mathbb{S}_2^1]]$--module. There is an exact sequence of complete $\Z_2[[\mathbb{S}_2^1]]$--modules
\begin{align*}
0 \ra \sC_3  \xra{\partial_3} \sC_2  \xra{\partial_2}\sC_1  \xra{\partial_1}\sC_0 \xra{\varepsilon} \Z_2 \ra 0, 
\end{align*}
where $\sC_0 \cong  \Z_2[[\mathbb{S}_2^1/G_{24}]]$ and $\sC_1 \cong  \sC_2 \cong  \Z_2[[\mathbb{S}_2^1/C_6]]$ and $\sC_3 = \Z_2[[\mathbb{S}_2^1/G_{24}']]$. Further, this is a free resolution of the trivial $\Z_2[[K^1]]$--module $\Z_2$.
\end{thm}

\begin{proof}
Let
\[\sC_3 :=\Z_2[[\mathbb{S}_2^1 / G_{24}']].\]
Let $\phi \co \sC_3 \ra N_2$ be the isomorphism of \fullref{thm:C3}.
Let $\partial_3 \co \sC_3 \ra \sC_2$ be the isomorphism $\phi$ followed by the inclusion of $N_2$ in $\sC_2$. This gives an exact sequence
\begin{equation}\label{resC2}
0 \ra \sC_3 \ra \sC_2 \xra{\partial_2} N_1 \ra 0.
\end{equation}
Splicing the exact sequences of (\ref{resN0}), (\ref{resN1}) and (\ref{resC2}) finishes the proof. 
\end{proof}

The exact sequence of \fullref{thm:res} is called the \emph{algebraic duality resolution}. The duality properties it satisfies will be described in \fullref{subsec:duality}.

\subsection{The algebraic duality resolution spectral sequence}

The algebraic duality resolution gives rise to a spectral sequence called \emph{algebraic duality resolution spectral sequence}, which we describe here. The following result is a refinement of \fullref{thm:DRSS}, which was stated in \fullref{subsec:state}. We define
\[Q_8' := \pi Q_8 \pi^{-1}.\]
We also let $V$ be the $2$--Sylow subgroup of $A_4$ and $V'$ be the $2$--Sylow subgroup of $A_4'$, where $A_4 \cong PG_{24}$ and $A_4'= PG_{24}'$ as defined in (\ref{A4}) and (\ref{A4'}).

\begin{thm}
Let $M$ be a profinite $\Z_2[[\mathbb{S}_2^1]]$--module. There is a first quadrant spectral sequence,
 \[E_1^{p,q} = \Ext^q_{\Z_2[[\mathbb{S}_2^1]]}(\sC_p, M) \Longrightarrow H^{p+q}(\mathbb{S}_2^1,M)\]
with differentials $d_r \co E_r^{p,q}\ra E_r^{p+r,q-r+1}$. Further:
 \begin{align*}
  E_1^{p,q}  \cong \begin{cases}
        H^q(G_{24}, M)  & \mbox{if $p= 0$}\\
    H^q(C_6,M)   & \mbox{if $p=1,2$}\\
     H^q(G_{24}', M)  & \mbox{if $p= 3$}\end{cases} 
    \end{align*}
Similarly, there are first quadrant spectral sequences
 \[E_1^{p,q} = \Ext^q_{\Z_2[[G]]}(\sC_p, M) \Longrightarrow H^{p+q}(G,M),\]
where $G$ is $S_2^1$, $P\mathbb{S}_2^1$ or $PS_2^1$. The $E_1$--term is
\[  E_1^{p,q}  \cong \begin{cases}
        H^p(Q_8; M)  & \mbox{if $q= 0$}\\
    H^p(C_2;M)   & \mbox{if $q=1,2$}\\
    H^p(Q_8'; M)  & \mbox{if $q= 3$}\end{cases}  \]
    when $G$ is $S_2^1$,
\[  E_1^{p,q}  \cong \begin{cases}
        H^p(A_4; M)  & \mbox{if $q= 0$}\\
    H^p(C_3;M)   & \mbox{if $q=1,2$}  \\
  H^p(A_4'; M)  & \mbox{if $q= 3$} \end{cases}   \] 
    when $G$ is $P\mathbb{S}_2^1$ and 
    \[  E_1^{p,q}  \cong \begin{cases}
        H^p(V; M)  & \mbox{if $q= 0$}\\
    H^p(\{e\};M)   & \mbox{if $q=1,2$}\\
    H^p(V'; M)  & \mbox{if $q= 3$} \end{cases}  \] 
    when $G$ is $PS_2^1$.
\end{thm}

\begin{proof}
There are two equivalent constructions. First, recall that the algebraic duality resolution
is spliced from the exact sequences
\begin{align} \label{N}
0 \ra N_i \ra \sC_i \ra N_{i-1} \ra 0,\end{align}
with $\sC_3=N_2$ and $N_{-1} = \Z_2$. The exact couple
\[\xymatrix{ \Ext(N_*, M) \ar[rr]^{\delta_*} & &  \Ext(N_{*-1}, M) \ar[dl]^{r_*}  \\ &  \Ext(\sC_*, M) \ar[ul]^{i^*} &   }\]
gives rise to the algebraic duality resolution spectral sequence. 

Alternatively, one can resolve each $\sC_{\bullet} \ra \Z_2$ into a double complex of projective finitely generated $\Z_2[[\mathbb{S}_2^1]]$--modules. The total complex $\Tot(P_{p,q})$ for $p \geq 0$ is a projective resolution of $\Z_2$ as a $\Z_2[[\mathbb{S}_2^1]]$--module. The homology of the double complex $\Hom_{\Z_2[[\mathbb{S}_2^1]]}( \Tot(P_{p,q}), M)$ is 
\[\Ext^{p+q}_{\Z_2[[\mathbb{S}_2^1]]}(\Z_2,M) \cong H^{p+q}(\mathbb{S}_2^1,M).\] The identification of the $E_1$--term follows from Shapiro's Lemma, which is \fullref{lem:shapiro} of the appendix. Indeed, any finite subgroup $H$ of $\mathbb{S}_2^1$ is closed. Further, since $\mathbb{S}_2^1 \cong S_2^1 \rtimes C_3$,
\begin{align*}\Ext^q_{\Z_2[[\mathbb{S}_2^1]]}(\Z_2[[\mathbb{S}_2^1]]\otimes_{\Z_2[H]}\Z_2,M) &\cong
\left(\Ext^q_{\Z_2[[{S}_2^1]]}(\Z_2[[{S}_2^1]]\otimes_{\Z_2[Syl_2(H)]}\Z_2,M)\right)^{C_3} \\
&\cong \left(\Ext^q_{\Z_2[Syl_2(H)]}(\Z_2,M)\right)^{C_3}  \\
&\cong H^q(H, M).
\end{align*}

For the groups $S_2^1$, $P\mathbb{S}_2$ and $PS_2^1$, one applies the same construction, keeping the following isomorphisms in mind. Let $H \subseteq \mathbb{S}_2^1$ be a finite subgroup which contains $C_6$ and let $PH =H/C_2$. There are isomorphisms
\begin{gather*}
\mathbb{Z}_2[[\mathbb{S}_2^1/H]] \cong \Z_2[[P\mathbb{S}_2^1]]\otimes_{\Z_2[PH]}\Z_2 \cong \Z_2[[P\mathbb{S}_2^1/PH]]  \\
\tag*{\text{and}}
\mathbb{Z}_2[[\mathbb{S}_2^1/H]] \cong \Z_2[[PS_2^1]]\otimes_{\Z_2[Syl_2(PH)]}\Z_2 \cong \Z_2[[P\mathbb{S}_2^1/Syl_2(PH)]]  \end{gather*}
as $\Z_2[[P\mathbb{S}_2^1]]$ and $\Z_2[[PS_2^1]]$--modules respectively.
\end{proof}

\subsection{The duality}\label{subsec:duality}
The algebraic duality resolution of \fullref{thm:res} owes its name to the fact that it satisfies a certain twisted duality. This duality is crucial for computations as it allows us to understand the map $\partial_3 \co \sC_3 \ra\sC_2$.

Let $\Mod({\mathbb{S}_2^1})$ denote the category of finitely generated left $\Z_2[[\mathbb{S}_2^1]]$--modules. Let $\pi=1+2\omega$ in $\mathbb{S}_2$ be as defined in (\ref{pi}). For $M$ in $\Mod(\mathbb{S}_2^1)$, let $c_{\pi}(M)$ denote the left $\Z_2[[\mathbb{S}_2^1]]$--module whose underlying $\Z_2$--module is $M$, but for which the action of $\gamma$ in $\mathbb{S}_2^1$ on an element $m$ in $c_{\pi}(M)$ is given by
\[\gamma \cdot m = \pi \gamma \pi^{-1} m.\]
If $\phi\co M \ra N$ is a morphism of left $\Z_2[[\mathbb{S}_2^1]]$--modules, let $c_{\pi}(\phi)\co c_{\pi}(M) \ra c_{\pi}(N)$ be given by
\[c_{\pi}(\phi)(m) = \phi(m).\]
Then $c_{\pi}\co \Mod({\mathbb{S}_2^1} )\ra \Mod({\mathbb{S}_2^1})$ is a functor.
In fact, $c_{\pi}$ is an involution, since $\pi^2 = -3$ is in the center of $\mathbb{S}_2$. 
We can now prove \fullref{thm:dualityintro}, which is restated here for convenience.

\begin{thm}[Henn--Karamanov--Mahowald, unpublished]\label{thm:resdual}
There exists an isomorphism of complexes of left $\Z_2[[\mathbb{S}_2^1]]$--modules:
\begin{align*}
\xymatrix{\ar @{} [dr]  0 \ar[r] &  \sC_3 \ar[d]^{f_3} \ar[r]^-{\partial_3} & \sC_2 \ar[d]^{f_2} \ar[r]^-{\partial_2} & \sC_1  \ar[d]^{f_1} \ar[r]^-{\partial_1} & \sC_0 \ar[d]^{f_0} \ar[r]^-{\varepsilon} & \mathbb{Z}_2 \ar@{=}[d] \ar[r] & 0  \\ 0 \ar[r]	&  c_{\pi}(\sC_0^*) \ar[r]^-{c_{\pi}(\partial_1^*)}	& c_{\pi}(\sC_1^*) \ar[r]^-{c_{\pi}(\partial_2^*)}	& c_{\pi}(\sC_2^*) \ar[r]^-{c_{\pi}(\partial_3^*)}	& c_{\pi}(\sC_3^*) \ar[r]^-{\overline{\varepsilon}}	& \mathbb{Z}_2 \ar[r]	 & 0	}
\end{align*}
\end{thm}
\begin{proof}
The proof is similar to the proof of Henn, Karamanov, Mahowald \cite[Proposition 3.8]{HKM}. Let $\sC_p^*= \Hom_{\Z_2[[\mathbb{S}_2^1]]}(\sC_p, \Z_2[[\mathbb{S}_2^1]])$ and $\partial_p^* = \Hom_{\Z_2[[\mathbb{S}_2^1]]}(\partial_p, \Z_2[[\mathbb{S}_2^1]])$ be the $\Z_2[[\mathbb{S}_2^1]]$--duals of $\sC_p$ and $\partial_p$ in the sense of \fullref{sec:contcoh}, (\ref{dualG}). The resolution of \fullref{thm:res} gives rise to a complex
\begin{align}\label{compK}
0 \ra \sC_0^* \xra{\partial_1^*} \sC_1^* \xra{\partial_2^*} \sC_2^* \xra{\partial_3^*} \sC_3^* \ra 0.
\end{align}
Because $K^1$ has finite index in $\mathbb{S}_2^1$, the induced and coinduced modules of $\Z_2[[K^1]]$ are isomorphic (see Symonds and Weigel \cite[Section 3.3]{SymWeig}). Therefore,
\[\Hom_{\Z_2[[\mathbb{S}_2^1]]}(\sC_p, \Z_2[[\mathbb{S}_2^1]]) \cong \Hom_{\Z_2[[K^1]]}(\sC_p, \Z_2[[K^1]]) \]
and the homology of the complex (\ref{compK}) is $H^n(K^1, \Z_2[[K^1]])$. By \fullref{cor:KP}, $H^n(K^1, \Z_2[[K^1]])$ is $0$ for $n\neq 3$ and $\Z_2$ for $n=3$. Further, the action of $G_{24}$ on $H^3(K^1, \Z_2[[K^1]])\cong \Z_2$ is trivial, as there are no non-trivial one dimensional $2$--adic representations of $G_{24}$. Hence, (\ref{compK}) is a resolution 
of $\Z_2$ as a trivial $\Z_2[[\mathbb{S}_2^1]]$--module.

The module $\sC_p^*$ is of the form $\Z_2[[\mathbb{S}_2^1/H]]$ via the isomorphism $t$ defined in \fullref{sec:contcoh}, (\ref{dual}). Let $\overline{\varepsilon}$ be the augmentation
\[\overline{\varepsilon} \co \sC_3^* \ra \Z_2. \]
Because the augmentation $\varepsilon\co \Z_2[[K^1]]\ra \Z_2$ induces an isomorphism 
\[\Hom_{\Z_2[[K^1]]}(\Z_2, \Z_2) \cong \Hom_{\Z_2[[K^1]]}(\Z_2[[K^1]], \Z_2),\]
one can choose an isomorphism $H^3(K, \Z_2[[K^1]]) \ra \Z_2$ making the following diagram commute:
\[\xymatrix{
\sC_3^* \ar[d]^-{\overline{\varepsilon}} \ar[r]
&H^3(K^1, \Z_2[[K^1]]) \ar@{.>}[ld]  \\
\Z_2}\]
Therefore, the dual resolution is given by
 \begin{align*}
0 \ra \sC_0^* \xra{\partial_1^*} \sC_1^* \xra{\partial_2^*} \sC_2^* \xra{\partial_3^*} \sC_3^* \xra{\overline{\varepsilon}} \Z_2 \ra 0.
\end{align*}
Take the image of this resolution in $\Mod(\mathbb{S}_2^1)$ under the involution $c_{\pi}$. Let $e_3^{\pi}$ be the canonical generator of $c_{\pi}(\sC_3^*)$. The map $f_0 \co \sC_0 \ra c_{\pi}(\sC_3^*)$ defined by
\[f_0(e_0) = e_3^{\pi}\]
is an isomorphism of $\Z_2[[\mathbb{S}_2^1]]$--modules, and the following diagram is commutative:
\[\xymatrix{\ar @{} [dr]   \sC_0    \ar[d]^{f_0}  \ar[r]^-{\varepsilon} & \Z_2 \ar@{=}[d]  \ar[r]  &0	  \\  c_{\pi}(\sC_3^*) \ar[r]^-{\overline{\varepsilon}} &\Z_2 \ar[r]  &0		}\]
Therefore, $f_0$ induces an isomorphism $\ker{\varepsilon} \cong \ker{\overline{\varepsilon}}$. As both
\begin{gather*}\sC_2 \xra{\ \ \ \ \partial_2\ \ \ \ } \sC_1 \xra{\ \ \ \ \partial_1 \ \ \ \ } \ker{\varepsilon}\\
\tag*{\text{and}}
c_{\pi}(\sC_2^*) \xra{\ c_{\pi}(\partial_2^*)\ } \sC_3^* \xra{\ c_{\pi}(\partial_3^*)\ } \ker{\overline{\varepsilon}}\end{gather*}
are the beginning of projective resolutions of $\ker{\varepsilon}$ and  $ \ker{\overline{\varepsilon}}$ as $\Z_2[[P\mathbb{S}_2^1]]$--modules, $f_0$ lifts to a chain map:
\[\xymatrix{\ar @{} [dr] 0 \ar[r] & \sC_3 \ar[d]^{f_3}    \ar[r]^-{\partial_3}   & \sC_2 \ar[d]^{f_2}    \ar[r]^-{\partial_2} & \sC_1 \ar[d]^{f_1} \ar[r]^-{\partial_1}  & \ker{\varepsilon}	\ar[d]^{f_0} \ar[r] & 0 \\  0 \ar[r] &  c_{\pi}(\sC_0^*) \ar[r]^-{c_{\pi}(\partial_1^*)}  & c_{\pi}(\sC_1^*) \ar[r]^-{c_{\pi}(\partial_2^*)}  & c_{\pi}(\sC_2^*)  \ar[r]^-{c_{\pi}(\partial_3^*)}  &  \ker{\overline{\varepsilon}}  \ar[r] & 0  } \]

Let $PS_2^1 = S_2^1/C_2$, where $S_2^1$ denotes the $2$--Sylow subgroup of $\mathbb{S}_2^1$. By construction, $f_0$ is an isomorphism, which implies that $\F_2 \otimes_{\Z_2[[PS_2^1]]} f_1$ and $\F_2 \otimes_{\Z_2[[PS_2^1]]} f_2$ are isomorphisms. As $\sC_p$ and $c_{\pi}(\sC_p^*)$ are projective $\Z_2[[P\mathbb{S}_2^1]]$--modules for $p=1,2$, \fullref{lem:ghmr} of the appendix implies that $f_1$ and $f_2$ are isomorphisms. Finally, $f_3$ must be an isomorphism by the five lemma.
\end{proof}

\subsection{A description of the maps}\label{subsec:maps}

This section is dedicated to proving the statements in \fullref{thm:maps}. The first statement of \fullref{thm:maps} is that
\[\partial_1(e_1) = (e-\alpha)e_0.\]
This was shown in \fullref{thm:res}. In this section, we prove the remaining statements of that theorem.

The following result provides a description of the maps $\partial_3 \co \sC_3 \ra \sC_2$ and proves the last part of \fullref{thm:maps}. It is a consequence of \fullref{thm:resdual}.
\begin{thm}\label{thm:d2}
There are isomorphisms of $\Z_2[[\mathbb{S}_2^1]]$--modules $g_p\co \sC_p  \ra \sC_p$ and differentials 
\[\partial_{p+1}' \co \sC_{p+1} \ra \sC_p,\]
such that
\begin{align}\label{res3}
\xymatrix{\ar @{} [dr]  0 \ar[r] &  \sC_3 \ar[d]^{g_3} \ar[r]^-{\partial_3} & \sC_2 \ar[d]^{g_2} \ar[r]^-{\partial_2} & \sC_1  \ar[d]^{g_1} \ar[r]^-{\partial_1} & \sC_0 \ar[d]^{g_0} \ar[r]^-{\varepsilon} & \mathbb{Z}_2 \ar@{=}[d] \ar[r] & 0  \\ 0 \ar[r]	&  \sC_3 \ar[r]^-{\partial_3'}	&  \sC_2 \ar[r]^-{\partial_2'}	&  \sC_1 \ar[r]^-{\partial_1'}		&\sC_0 \ar[r]^-{\varepsilon}	& \mathbb{Z}_2 \ar[r]	 & 0	}
\end{align}
is an isomorphism of complexes.
The map $\partial_3'\co \sC_3 \ra \sC_2$ is given by
\begin{align}\label{p2'}
\partial_3'(e_3) = \pi (e+i+j+k)(e-\alpha^{-1}) \pi^{-1} e_2.
\end{align}
\end{thm}

\begin{proof}
We will construct a commutative diagram:
\begin{align*}
\xymatrix{\ar @{} [dr]  0 \ar[r] &  \sC_3 \ar[d]^{f_3} \ar[r]^-{\partial_3} & \sC_2 \ar[d]^{f_2} \ar[r]^-{\partial_2} & \sC_1  \ar[d]^{f_1} \ar[r]^-{\partial_1} & \sC_0 \ar[d]^{f_0} \ar[r]^-{\varepsilon} & \mathbb{Z}_2 \ar@{=}[d] \ar[r] & 0  \\ 0 \ar[r]	&  c_{\pi}(\sC_0^*) \ar[d]^{q_3}  \ar[r]^-{c_{\pi}(\partial_1^*)}	& c_{\pi}(\sC_1^*) \ar[d]^{q_2} \ar[r]^-{c_{\pi}(\partial_2^*)}	& c_{\pi}(\sC_2^*) \ar[d]^{q_1}  \ar[r]^-{c_{\pi}(\partial_3^*)}	& c_{\pi}(\sC_3^*) \ar[d]^{q_0}  \ar[r]^-{\varepsilon}	& \mathbb{Z}_2 \ar@{=}[d] \ar[r]	 & 0\\
0 \ar[r]	&  \sC_3 \ar[r]^-{\partial_3'}	&  \sC_2 \ar[r]^-{\partial_2'}	&  \sC_1 \ar[r]^-{\partial_1'}		&\sC_0 \ar[r]^-{\varepsilon}	& \mathbb{Z}_2 \ar[r]	 & 0 }
\end{align*}
The maps $g_p$ will be the composites of the vertical maps. First, let $e_p^{\pi} \in c_{\pi}(\sC_p^*)$ be the canonical generator. Define isomorphisms $q_{p}\co c_{\pi}(M_{3-p}^*) \ra M_{p}$ by 
\[q_{p}(e_{3-p}^{\pi}) = e_p.\]
Define $g_p\co \sC_p \ra \sC_p$ by 
\[g_p = q_pf_p\] and $\partial_{p+1}' \co \sC_{p+1}\ra \sC_p$ by
\[\partial_{p+1}' = q_p c_{\pi}(\partial_{3-p}^*) q_{p+1}^{-1}.\]
By construction, (\ref{res3}) is commutative.

In order to compute $\partial_3'$, it is necessary to understand $\partial_1^*$. By definition,
\[ \partial_1^{*}(e_0^{*})(e_1)= e_0^{*}((e-\alpha)e_1) = (e-\alpha) \sum_{h \in G_{24}} h . \]
However,
\begin{align*}
(e-\alpha) \sum_{h \in G_{24}} h   &= (e-\alpha) \sum_{h \in C_{6}}h (e+i^{-1}+j^{-1}+k^{-1}) \\
&=  \sum_{h \in C_{6}} h (e-\alpha) (e+i^{-1}+j^{-1}+k^{-1})\\
&=\left((e+i+j+k)(e-\alpha^{-1})e_1^{*}\right)(e_1).
\end{align*}
Hence,
\[\partial_1^{*}(e_0^*) = (e+i+j+k)(e-\alpha^{-1}) e_1^*.\]
A diagram chase shows that $\partial_3'$ is given by (\ref{p2'}).
\end{proof}

The maps $\partial_1 \co \sC_1 \ra \sC_0$ and $\partial_3\co \sC_3 \ra \sC_2$ now have explicit descriptions up to isomorphisms. The map $\partial_2 \co  \sC_2 \ra \sC_1$ is harder to describe. \fullref{thm:BETAmod8} and \fullref{cor:beta} below give an estimate for this map. These are technical results which will be used in our computations in \cite{Paper2}. Note that \fullref{thm:res}, \fullref{thm:d2} and \fullref{cor:beta} below prove \fullref{thm:maps}, which was stated in \fullref{subsec:state}.

Recall that 
\[\alpha_{\tau} = [\tau, \alpha] = \tau \alpha \tau^{-1}\alpha^{-1}.\]
We will need the following result to describe the element $\Theta$ of \fullref{lem:N2}. 
\begin{lem}\label{lem:62}
Let $n\geq 2$ and $x$ be in $I F_{n/2} K^1$. There exist $h_0$, $h_1$ and $h_2$ in $\Z_2[[F_{n/2}K^1]]$ such that:
\begin{align}\label{xform}
x=\begin{cases}
h_0(e-\alpha^{2^{m-1}}) + h_1(e-\alpha_i^{2^{m-1}})+h_2(e-\alpha_j^{2^{m-1}}) & \mbox{if $n=2m$}\\
         h_0(e-\alpha^{2^{m}}) + h_1(e-\alpha_i^{2^{m-1}})+h_2(e-\alpha_j^{2^{m-1}})& \mbox{if $n=2m+1$}\end{cases}\end{align} 
\end{lem}
\begin{proof}
Define a map of $\Z_2[[F_{n/2}K^1]]$--modules
\begin{align*}
\displaystyle p \co \bigoplus_{i=0}^{2} \Z_2[[F_{n/2}K^1]]_i \ra I F_{n/2}K^1
\end{align*}
by sending $(h_0, h_1,h_2)$ to the element given by (\ref{xform}). It is sufficient to show that the map induced by $p$ surjects onto
\[H_1(F_{n/2}K^1, \F_2) \cong \F_2 \otimes_{ \Z_2[[F_{n/2}K^1]]}I F_{n/2}K^1.\]
By \fullref{lem:cent}, $H_1(F_{n/2}K^1,\F_2)$ is generated by the classes $\alpha^{2^{m-1}}$, $\alpha_i^{2^{m-1}}$ and $\alpha_j^{2^{m-1}}$ if $n=2m$, and the classes $\alpha^{2^{m}}$, $\alpha_i^{2^{m-1}}$ and $\alpha_j^{2^{m-1}}$ if $n=2m+1$. Therefore, $\F_2 \otimes_{ \Z_2[[F_{n/2}K^1]]}p$ is surjective, and hence so is $p$.
\end{proof}
The ideal
\begin{equation*}\mathcal{I} = ((IK^1)^7, 2(IK^1)^3, 4 (IK^1), 8)\end{equation*}
will play a crucial role in the following estimates. 
\begin{cor}\label{cor:finalcor}
Let $e_0$ be the canonical generator of $\sC_0$ and $g$ be in $F_{8/2}K^1$. There exists $h$ in $\Z_2[[\mathbb{S}_2^1]]$ such that
$(e-g)e_0 = h(e-\alpha)e_0$
with $h \equiv 0  \mod \mathcal{I}$. 
\end{cor}
\begin{proof}
By \fullref{lem:62}, there exist $h_0$, $h_1$ and $h_2$ in $\Z_2[[F_{8/2}K^1]]$ such that
 \[e-g= h_0(e-\alpha^{8}) + h_1(e-\alpha_i^{8})+h_2(e-\alpha_j^{8}).\]
Since
\[ (e-x^8) = \left(\sum_{s=0}^{7}x^s\right)(e-x) \]
this implies that
\begin{align*}
e-g &=
h_0\left(\sum_{s=0}^{7}\alpha^s\right)(e-\alpha) +h_1 \left(\sum_{s=0}^{7}\alpha_i^s\right)(e-\alpha_i) +h_2\left(\sum_{s=0}^{7}\alpha_j^s\right)(e-\alpha_j) .
\end{align*}
Let
\[h =
h_0\left(\sum_{s=0}^{7}\alpha^s\right) + h_1 \left(\sum_{s=0}^{7}\alpha_i^s\right)(i-\alpha_i) +h_2 \left(\sum_{s=0}^{7}\alpha_j^s\right)(j-\alpha_j).\]
If $\tau \in G_{24}$, then $\tau e_0 = e_0$. Hence,
\[(\tau-\alpha_{\tau})(e-\alpha)e_0 =(e-\alpha_\tau)e_0.\]
Using this fact, one verifies that $(e-g) e_0 = h(e-\alpha)e_0$. Further,
\[\sum_{s=0}^{7}x^s \equiv (1-x)^7+2x^4(x-1)^3+4x^2(x-1)  \mod (8). \]
Since $\alpha$, $\alpha_i$ and $\alpha_j$ are in $K^1$ and $K^1$ is a normal subgroup, this implies that 
\[h \equiv 0  \mod ((IK^1)^7, 2(IK^1)^3, 4 (IK^1), 8). \proved \]
\end{proof}

We will use the following result.
\begin{lem}\label{lem:deepfilt}
The element $\alpha_i\alpha_j\alpha_k$ is in $F_{4/2}K^1$. The element $\alpha_i\alpha_j\alpha_k\alpha^2$ is in $F_{8/2}K^1$.
\end{lem}
\begin{proof}
Let $T = \alpha S$ in $\mathcal{O}_2 \cong \End(F_2)$. Then $T^2 = -2$, and $a T = T a^{\sigma}$ for $a$ in $\W$. As defined in (\ref{alpha}) and \fullref{lem:Q8split}, we have:
\begin{align*}
\alpha &= \frac{1}{\sqrt{-7}} (1-2\omega) & i &= -\frac{1}{3}(1+2\omega)(1-T) \\ 
 j &= -\frac{1}{3}(1+2\omega)(1-\omega^2T) & k&=  -\frac{1}{3}(1+2\omega)(1-\omega T)
\end{align*}
Further, 
\begin{align*}
\alpha^{-1} &= -\frac{1}{\sqrt{-7}} (1-2\omega^2).
\end{align*}
We use the fact that $\frac{1}{3}$ and $ \frac{1}{\sqrt{-7}}$ are in $Z(\mathbb{S}_2)$. We also use the fact $\tau^{-1}=-\tau$ for $\tau=i$, $j$ and $k$ and the fact that $S^4 = 4$ and $S^8 = 16$.

First, note that
\begin{align*}
i\alpha &= - \frac{1}{3\sqrt{-7}}(1+2\omega)(1-T)   (1-2\omega) \\
&= - \frac{1}{3\sqrt{-7}}(1+2\omega) (  (1-2\omega)-(1-2\omega^2)T) \\
&= - \frac{1}{3\sqrt{-7}}((5+4\omega)+(1-4\omega)T).
\end{align*}
Further,
\begin{align*}
i^{-1}\alpha^{-1} &=  -\frac{1}{3\sqrt{-7}}(1+2\omega)(1-T)  (1-2\omega^2) \\
&=  -\frac{1}{3\sqrt{-7}}(1+2\omega)((1-2\omega^2)-(1-2\omega)T)  \\
&=  -\frac{1}{3\sqrt{-7}}((-1+4\omega)-(5+4\omega)T) .
\end{align*}
Therefore,
\begin{align*}
\alpha_i &= i \alpha i^{-1}\alpha^{-1} \\
&=  -\frac{1}{63}((5+4\omega)+(1-4\omega)T)((-1+4\omega)-(5+4\omega)T) \\
&\equiv 13+(2+8\omega)T \mod S^8. 
\end{align*}
Using the fact that $\alpha_{j} = \omega \alpha_i \omega^2$ and $\alpha_{k} = \omega^2 \alpha_i \omega$, this implies that:
\begin{align*}
\alpha_j &\equiv 13+\omega^2(2+8\omega)T \mod S^8 &
\alpha_k &\equiv 13+\omega (2+8\omega)T \mod S^8
\end{align*}
Hence,
\begin{align*}
\alpha_i\alpha_j &\equiv (13+(2+8\omega)T) (13+\omega^2 (2+8\omega)T )\\
&\equiv \big(9+\omega^2 (10+8\omega)T+(10+8\omega)T + (2+8\omega)(\omega (2+8\omega^2))T^2\big)\\
&\equiv 9+8\omega+(8+14\omega)T \mod S^8,
\end{align*}
so that
\begin{align*}
\alpha_i\alpha_j\alpha_k &\equiv (9+8\omega+(8+14\omega)T )(13+\omega(2+8\omega)T) \\
&\equiv (5+8\omega+(8+6\omega)T+(9+8\omega)\omega (2+8\omega)T + (8+14\omega)\omega^2 (2+8\omega^2)T^2 ) \\
&\equiv 13+8\omega \mod S^8.
\end{align*}
This shows that $\alpha_i\alpha_j\alpha_k \equiv 1$ modulo $S^4$. Finally, note that
\begin{align*}
\alpha_i\alpha_j\alpha_k\alpha^2 &\equiv (13+8\omega)\left(\frac{1}{\sqrt{-7}} (1-2\omega)\right)^2 \\
&\equiv -\frac{1}{7}(13+8\omega)(1-2\omega)^2 \\
&\equiv -\frac{9}{7} \\
&\equiv 1 \mod S^8,
\end{align*}
which shows that $\alpha_i\alpha_j\alpha_k\alpha^2$ is in $F_{8/2}K^1$.
\end{proof}

\begin{thm}\label{thm:BETAmod8}
There is an element $\Theta$ in $\Z_2[[\mathbb{S}_2^1]]$ satisfying the conditions of \fullref{lem:N2} such that 
\begin{align*}
\Theta  &\equiv e + \alpha+i+j+k -\alpha_i-\alpha_j -\alpha_k \\
& \ \ \ -\frac{1}{3} tr_{C_3}\bigg((e-\alpha_i)(j-\alpha_j)+(e-\alpha_i\alpha_j) (k-\alpha_k)+(e-\alpha_i\alpha_j\alpha_k)(e+\alpha)\bigg)
\end{align*}
modulo $\mathcal{I}= ((IK^1)^7, 2(IK^1)^3, 4 (IK^1), 8)$, where $ tr_{C_3}$ is defined by (\ref{trace}).
\end{thm}

\begin{proof}
We will use the following facts. First, note that
\[\tau e_q = e_q\]
for $\tau \in G_{24}$ and $q=0$, or for $\tau \in C_6$ and $q=1$. This implies that
\[\tau (e-\alpha) e_0 =(e- \alpha_{\tau} \alpha) e_0. \]
Since $j=\omega i \omega^{-1}$ and $k = \omega^{-1} i\omega$, it also implies:
\begin{align*}\omega i e_q &= j e_q &  \omega^2i e_q &= ke_q\end{align*}
The element $\alpha \in \W^{\times}\subseteq \mathbb{S}_2$ commutes with $\omega$. This implies that 
\[\omega \alpha_{i}e_q= \alpha_j e_q.\] 
We will use the fact that for $\tau \in G_{24}$,
\begin{align*}
(\tau-\alpha_{\tau})(e-\alpha)e_0 =(e-\alpha_\tau)e_0.
\end{align*}
We will also use the following identity:
\begin{align*}
e-gh = (e-g)+(e-h)-(e-g)(e-h)\end{align*}

Let $\Theta_0 = e+i$. Then $tr_{C_3}(\Theta_0)e_1 = (3+i+j+k)e_1$ and
\begin{align*}
\partial_{1}(\Theta_0e_1)&= ( e+i)(e-\alpha)e_0 \\
&= (e-\alpha)e_0 + (e-\alpha_i\alpha)e_0 \\
&=2 (e-\alpha)e_0 + (e-\alpha_i)e_0 - (e-\alpha_i)(e-\alpha)e_0 \\
&=(e-\alpha^2)e_0+(e-\alpha)^2e_0  + (e-\alpha_i)e_0- (e-\alpha_i)(e-\alpha)e_0 .
\end{align*}
Let $\Theta_1 = e+i -(e-\alpha)+(e-\alpha_i)$. Then
\begin{align*}
\partial_{1}(\Theta_1e_1)&= (e-\alpha^2)e_0  + (e-\alpha_i)e_0.
\end{align*}
Therefore,
\begin{align*}
\partial_{1}(tr_{C_3}(\Theta_1)e_1)&=3(e-\alpha^2)e_0  + (e-\alpha_i)e_0+ (e-\alpha_j)e_0+ (e-\alpha_k)e_0 \\
&= 3(e-\alpha^2)e_0+ (e-\alpha_i)(e-\alpha_j)e_0+(e-\alpha_i\alpha_j)e_0+ (e-\alpha_k)e_0 \\
&= 3(e-\alpha^2)e_0+ (e-\alpha_i)(e-\alpha_j)e_0+(e-\alpha_i\alpha_j) (e-\alpha_k)e_0\\
& \ \ \ \ \ +(e-\alpha_i\alpha_j\alpha_k)e_0  \\ 
&= 2(e-\alpha^2)e_0+ (e-\alpha_i)(e-\alpha_j)e_0+(e-\alpha_i\alpha_j) (e-\alpha_k)e_0\\
&\ \ \ \ \ +
(e-\alpha_i\alpha_j\alpha_k)(e-\alpha^2)e_0+
(e-\alpha_i\alpha_j\alpha_k\alpha^2)e_0  .
\end{align*}
Let
\begin{align*}
\Theta_2 &= tr_{C_3}\big(e+i -(e-\alpha)+(e-\alpha_i) \big) -2(e+\alpha) \\
& \ \ \ \ \ -(e-\alpha_i)(j-\alpha_j)-(e-\alpha_i\alpha_j) (k-\alpha_k)-(e-\alpha_i\alpha_j\alpha_k)(e+\alpha).
\end{align*}
Then $\Theta_2 \equiv 3+i+j+k \mod (4, IK^1)$. Further,
\begin{align*}
\partial_1(\Theta_2e_1) =(e-\alpha_i\alpha_j\alpha_k\alpha^2)e_0.
\end{align*}
By \fullref{lem:deepfilt}, $\alpha_i\alpha_j\alpha_k\alpha^2 \in F_{8/2}K^1$. By \fullref{cor:finalcor}, there exists $h$ such that
\[(e-\alpha_i\alpha_j\alpha_k\alpha^2)e_0 =h(e-\alpha)e_0\]
and $h \equiv 0$ modulo $\mathcal{I}$, where  $\mathcal{I} = ((IK^1)^7, 2(IK^1)^3, 4 IK^1, 8)$. Therefore,
\[\partial_1((\Theta_2  -h)e_1) = (e-\alpha_i\alpha_j\alpha_k\alpha^2)e_0 -h(e-\alpha)e_0 = 0.\]
Define
\begin{align*}
\Theta  =  \frac{1}{3} tr_{C_3}(\Theta_2-h).
\end{align*}
Then $\Theta$ satisfies the conditions of \fullref{lem:N2}. Further,
\begin{align*}
\Theta &\equiv e + \alpha+i+j+k -\alpha_i-\alpha_j -\alpha_k \\
& \ \ \ -\frac{1}{3} tr_{C_3}\bigg((e-\alpha_i)(j-\alpha_j)+(e-\alpha_i\alpha_j) (k-\alpha_k)+(e-\alpha_i\alpha_j\alpha_k)(e+\alpha)\bigg)
\end{align*}
modulo $\mathcal{I}$.
\end{proof}

\begin{cor}\label{cor:beta}
Let $\mathcal{J} = (IF_{4/2}K^1, (IF_{3/2}K^1)(IS_2^1), \mathcal{I})$. The element $\Theta$ of \fullref{thm:BETAmod8} satisfies
\[\Theta \equiv  e+\alpha+i+j+k-\alpha_i-\alpha_j-\alpha_k \mod \mathcal{J}\]
and $\Theta \equiv e+\alpha$ modulo $(2, (IS_2^1)^2)$.
\end{cor}
\begin{proof}
First, note that $\alpha_{\tau} \in F_{3/2}K^1$ for $\tau\in G_{24}$. Further, by \fullref{lem:deepfilt}, $\alpha_i\alpha_j\alpha_k$ is in $F_{4/2}K^1$. Hence, it follows from \fullref{thm:BETAmod8} that
\[ \Theta \equiv e+\alpha+i+j+k-\alpha_i-\alpha_j-\alpha_k \mod \mathcal{J}.\]
For the second claim, we first prove that $\mathcal{J} \subseteq (2, (IS_2^1)^2)$. It is clear that,
\[ ((IF_{3/2}K^1)(IS_2^1), \mathcal{I})\subseteq (2, (IS_2^1)^2).\]
Further, it follows from \fullref{lem:62} and the fact that $(e-x^{2^k}) \equiv (e-x)^{2^k} $ modulo $(2)$ that 
\[IF_{4/2}K^1 \subseteq (2, (IS_2^1)^2 ).\]
Therefore, $\mathcal{J} \subseteq (2,(IS_2^1)^2)$. Hence,
\begin{align*}
\Theta &\equiv e+\alpha+i+j+k-\alpha_i-\alpha_j-\alpha_k  \mod (2,  (IS_2^1)^2).
 \end{align*}
Further, $(e-i)(e-j)\equiv  e+i+j+k$ modulo $(2)$ and
\begin{align*}
e-\alpha_i   &= i\alpha((e-\alpha^{-1})(e-i^{-1})  -(e-i^{-1})(e-\alpha^{-1})).
\end{align*}
Therefore, $e+i+j+k$ and $e-\alpha_{\tau}$ are in $(2,(IS_2^1)^2)$. We conclude that
\[\Theta \equiv e+\alpha \mod (2,  (IS_2^1)^2) . \proved \]
\end{proof}

 \section{Appendix: Background on profinite groups}\label{sec:contcoh}
 
We use the terminology of Ribes and Zalesskii \cite[Section 5]{zalesskii}. Let $G$ be a profinite $p$--adic analytic group and $\{U_k\}$ be a system of open normal subgroups of $G$ such that $\bigcap_k U_k =\{e\}$. The \emph{completed group ring of $G$} is
\begin{align*}
\Z_p[[G]] := \displaystyle \lim_{n,k} \Z/(p^n)[G/U_k]. 
\end{align*}
The \emph{augmentation} is the continuous homomorphism of $\Z_p$--modules $\varepsilon \co \Z_p[[G]]  \ra \Z_p$ defined by $\varepsilon(g)=1$ for $g \in G$. The \emph{augmentation ideal} $IG$ is the kernel of $\varepsilon$. 

A \emph{left $\Z_p[[G]]$--module} is a $\Z_p[[G]]$--module $M$ which is a Hausdorff topological abelian group with a continuous structure map $\Z_p[[G]] \times M \ra M$. The module $M$ is \emph{finitely generated} if it is the closure of the $\Z_p[[G]]$--module generated by a finite subset of $M$. It is \emph{discrete} if it is the union of its finite $\Z_p[[G]]$--submodules and \emph{profinite} if it is the inverse limit of its finite $\Z_p[[G]]$--submodule quotients (\cite[Lemma 5.1.1]{zalesskii}). 
The module $M$ is complete with respect to the $IG$--adic topology if 
\[M \cong  \displaystyle \lim_{n,k} \Z_p/(p^n)[G/U_k] \otimes_{\Z_p[[G]]}M.\]
It is a theorem of Lazard that $\Z_p[[G]]$ is Noetherian (see Symonds and Weigel \cite[Theorem 5.1.2]{SymWeig}). Finitely generated $\Z_p[[G]]$--modules are thus both profinite and complete with respect to the $IG$--adic topology. 

Let $M = \lim_i M_i$ be a profinite $\Z_p[[G]]$--bimodule and $N =\lim_j N_j$ a profinite left $\Z_p[[G]]$--module. Then
\[M \otimes_{\Z_p[[G]]}N = \lim_{i,j} M_i \otimes_{\Z_p[[G]]} N_j\] 
denotes the completed tensor product, which is itself a profinite left $\Z_p[[G]]$--module \cite[Section 5.5]{zalesskii}.
The abelian group of continuous $\Z_p[[G]]$--homomorphisms is denoted by 
\[\Hom_{\Z_p[[G]]}(M, N).\] 
This is a topological space with the compact open topology. If $M$ is finitely generated, then it is compact \cite[Section 3.7]{SymWeig}.

Lazard also proves in \cite[V.3.2.7]{lazard} that the trivial $\Z_p[[G]]$--module $\Z_p$ admits a resolution by finitely generated $\Z_p[[G]]$--modules. A $\Z_p[[G]]$--module $M$ which admits a projective resolution $P_{\bullet} \ra M$ by finitely generated $\Z_p[[G]]$--modules is said to be of type ${\bf FP^{\infty}}$ \cite[Section 3.7]{SymWeig}. For such $M$, we let
\begin{gather*}\Ext_{\Z_p[[G]]}^n (M, N) = H^n(\Hom_{\Z_p[[G]]}(P_{\bullet}, N))\\
\tag*{\text{and}}
\Tor^{\Z_p[[G]]}_n(M, N) = H_n( P_{\bullet} \otimes_{\Z_p[[G]]} N).\end{gather*}
These functors are studied by Symonds and Weigel in \cite[Section 3.7]{SymWeig}. There are isomorphisms
\[H^n(G, N)   \cong \Ext_{\Z_p[[G]]}^n (M, N),  \]
where $H^n(G, N)$ is the cohomology computed with continuous cochains and
\[H_n(G, N)   \cong \Tor^{\Z_p[[G]]}_n(\Z_p, N) \]
(see Neukirch, Schmidt and Wingberg \cite[Propositions 5.2.6, 5.2.14]{neukirch} or the discussion in Kohlhaase \cite[Section 3]{kohlhaase}). Therefore, these functors satisfy the usual properties of group cohomology (see Ribes and Zalesskii \cite[Section 6]{zalesskii}). In particular, for $[G,G]$ the commutator subgroup, $G^* = G^p[G,G]$, and $\Z_p$ and $\F_p$ the trivial modules, we have:
\begin{align*}
H_1(G, \Z_p) &\cong  G/\overline{[G,G]} &  H_1(G, \F_p) &\cong G/\overline{G^*} \\
H^1(G, \Z_p) &\cong \Hom(G, \Z_p)&  H^1(G, \F_p) &\cong \Hom(G, \F_p) \end{align*}

\begin{exmps}We give examples which we use in this paper, in \cite{Paper2} and in \cite{Paper5}.
\begin{enumerate}[(a)]
\item The modules 
\[\Z_p[[G/H]]:= \Z_p[[G]] \otimes_{\Z_p[H]}\Z_p\]
for $H$ a finite subgroup of $G$ and $\Z_p$ the trivial $\Z_p[H]$--module are finitely generated, and thus profinite and complete. 
\item The $\Z_p[[G]]$--dual of a finitely generated $\Z_p[[G]]$--module $M$ is defined as
\begin{align}\label{dualG}
M^* := \Hom_{\Z_p[[G]]}(M, \Z_p[[G]])
\end{align}
with the action of $g \in G$ on $\phi \in M^*$ defined by
\[(g\phi)(m)=\phi(m)g^{-1}.\]
This gives $M^*$ the structure of a finitely generated left $\Z_p[[G]]$--module (see Symonds and Weigel \cite[3.7.1]{SymWeig} and Henn, Karamanov and Mahowald \cite[Section 3.4]{HKM}). 
For example, if $H \subseteq G$ is a finite subgroup and $[g]$ denotes the coset $gH$, there is a canonical isomorphism 
\begin{align}\label{dual}
t\co \Z_p[[G/H]] \ra \Z_p[[G/H]]^*
\end{align}
which sends $[g]$ to the map $[g]^*\co \Z_p[[G/H]] \ra \Z_p[[G]]$ defined by
\[[g]^*([x]) = x \sum_{h \in H} h g^{-1}.\]
We refer the reader to \cite[Section 3.4]{HKM} for a detailed discussion of $\Z_p[[G]]$--duals.
\item In the case when $G=S_n$ is the $p$--Sylow subgroup of $\mathbb{S}_n$, an important example is the continuous $\Z_p[[S_n]]$--module $(E_n)_*X =\pi_*L_{K(n)}(E_n\smsh X)$ for a spectrum $X$ (see Goerss, Henn, Mahowald and Rezk \cite[Section 2]{ghmr}). In the case when $X$ is a finite spectrum, $(E_n)_*X$ is profinite, although it is not known if, in general, it is finitely generated over $\Z_p[[S_n]]$. For a more extensive discussion, see the work of Kohlhaase in \cite{kohlhaase}.
\end{enumerate}
\end{exmps}

\begin{lem}[Shapiro's Lemma]\label{lem:shapiro}
Let $G$ be a profinite $p$--analytic group and let $H$ be a closed subgroup. Let $M$ be a $\Z_p[[H]]$--module of type ${\bf FP^{\infty}}$ and let $N = \lim_i N_i$ be a profinite $\Z_p[[G]]$--module, which is also a $\Z_p[[H]]$--module via restriction. Then,
\[\Ext_{\Z_p[[G]]}^*(\Z_p[[G]] \otimes_{\Z_p[[H]]} M, N) \cong \Ext_{\Z_p[[H]]}^*(M, N) . \]
\end{lem}
\begin{proof}
Let $P_{\bullet}\ra M$ be a projective resolution of $M$ by finitely generated $\Z_p[[H]]$--modules. According to Brumer \cite[Lemma 4.5]{brumer}, $\Z_p[[G]]$ is a projective $\Z_p[[H]]$--module. Hence, the functor $\Z_p[[G]]\otimes_{\Z_p[[H]]}(-)$ is exact, and $\Z_p[[G]]\otimes_{\Z_p[[H]]} P_{\bullet}$ is a projective resolution of $\Z_p[[G]]\otimes_{\Z_p[[H]]} M$ by finitely generated $\Z_p[[G]]$--modules. Finally, note that
\begin{align*}
\Hom_{\Z_p[[G]]}(\Z_p[[G]]\otimes_{\Z_p[[H]]} P_{\bullet}, N) &\cong \displaystyle \lim_i \Hom_{\Z_p[[G]]}(\Z_p[[G]]\otimes_{\Z_p[[H]]} P_{\bullet},  N_i) \\
&\cong \displaystyle \lim_i \Hom_{\Z_p[[H]]}( P_{\bullet},  N_i) \\
 &\cong  \Hom_{\Z_p[[H]]}( P_{\bullet},  N),
\end{align*}
where the first isomorphism is proved by Symonds and Weigel \cite[(3.7.1)]{SymWeig} and the second follows from Ribes and Zalesskii \cite[Proposition 5.5.4 (c)]{zalesskii}. \end{proof}

The following result is Lemma 4.3 of Goerss, Henn, Mahowald and Rezk \cite{ghmr}. It is a version of Nakayama's lemma in this setting.
\begin{lem}[Goerss--Henn--Mahowald--Rezk]\label{lem:ghmr}
Let $G$ be a finitely generated profinite $p$--group. Let $M$ and $N$ be finitely generated complete $\Z_p[[G]]$--modules and $f \co M\rightarrow N$ be a map of complete $\Z_p[[G]]$--modules. If the induced map 
\[\F_p\otimes_{\Z_p[[G]]} f \co \F_p \otimes_{\Z_p[[G]]} M \ra \F_p \otimes_{\Z_p[[G]]} N\]
is surjective, then so is $f$. If the map
\[\Tor_{q}^{\Z_p[[G]]}(\F_p, f)  \co \Tor_{q}^{\Z_p[[G]]}(\F_p, M) \ra \Tor_{q}^{\Z_p[[G]]}(\F_p, N)\]
is an isomorphism for $q=0$ and surjective for $q=1$, then $f$ is an isomorphism.
\end{lem}

The following is a restatement of some of the results which can be found in Ribes and Zalesskii \cite[Lemma 6.8.6]{zalesskii}
\begin{lem}\label{lem:Ggen}
Let $G$ be a profinite group and let $IG$ be the augmentation ideal. For a profinite $\Z_p[[G]]$--module $M$, the boundary map for the short exact sequence
\[0 \ra IG \ra \Z_p[[G]] \xra{\varepsilon} \Z_p \ra 0\]
induces an isomorphism
\[H_{n+1}(G,M) \cong \Tor_n^{\Z_p[[G]]}(IG,M).\]
For the trivial module $M = \Z_p$, this isomorphism sends $g$ in $G/\overline{[G,G]}$ to the residue class of $e-g$ in $H_0(G, IG)\cong IG/\overline{IG^2}$. Let $G^*$ be the subgroup generated by $[G,G]$ and $G^p$. For $M=\F_p$, it sends $g$ in  $G/\overline{G^*}$ to the residue class of $e-g$ in $\F_p \otimes_{\Z_p} IG/\overline{IG^2}$
\end{lem}

Finally, we note the following classical result.
\begin{lem}\label{lem:sq1}
Let $G$ be a profinite $2$--analytic group and suppose that $H_1(G,\Z_2) \cong G/\overline{[G,G]}$ is a finitely generated $2$--group. Suppose that the residue class of an element $g$ in $G/\overline{[G,G]}$ generates a summand isomorphic to $\Z/2^k$. Let $x$ in $H^1(G ,\Z/2) \cong \Hom(G, \Z/2)$ be the homomorphism dual to $g$. Then $x^2$ is non-zero in $H^2(G, \Z/2)$ if and only if $k=1$. 
\end{lem}
\begin{proof}
This follows from the fact that $x$ in $H^1(G, \Z/2)$ has a non-zero Bockstein in the long exact sequence associated to the extension of trivial modules 
\[1 \ra \Z/2 \ra \Z/4 \ra \Z/2 \ra 1\] 
if and only if $g$ generates a $\Z/2$ summand. 
\end{proof}

%%%%%%%%%%%%%%%%%%%%   End of main body of article
%
%                             References
%
%   BiBTeX users uncomment the following line:
%
%\bibliographystyle{gtart}
%

\end{document}